\documentclass{amsart}
\usepackage{amssymb}
\usepackage{amsmath}
\usepackage{amsfonts}

\newtheorem{theorem}{Theorem}
\theoremstyle{plain}

\newtheorem{corollary}{Corollary}

\newtheorem{definition}{Definition}

\newtheorem{lemma}{Lemma}
\newtheorem{notation}{Notation}

\numberwithin{equation}{section}
\input{tcilatex}

\begin{document}
\title{On a new class of systems of generalized quasy-variational
inequalities }
\author{Monica Patriche}
\maketitle

Abstract: We introduce new types of systems of generalized quasi-variational
inequalities and we prove the existence of the solutions by using results of
pair equilibrium existence for free abstract economies. We consider the
fuzzy models and we also introduce the random free abstract economy and the
random equilibrium pair. The existence of the solutions for the systems of
quasi-variational inequalities comes as consequences of the existence of
equilibrium pairs for the considered free abstract economies.


\section{Introduction}

Variational inequality theory proved to be a powerful tool used to formulate
a variety of equilibrium problems concerning the traffic network, the
spatial price, the oligopolistic market, the financial domain or the
migration. It also concerns the qualitative analysis of the problems in
terms of existence and uniqueness of solutions, stability and sensitivity
analysis and it provides algorithms for computational purposes. Systems of
non-linear equations, optimization problems, complementarity problems, fixed
point theorems are contained as special cases of variational inequality
theory. Hartman and Stampacchia (1966) introduced this domain as a tool for
the study of partial differential equations with applications principally
drawn from mechanics. Since then, the theory of variational inequalities has
been developing very fast. New results have been obtained, for instance, in
[2]-[4], [6], [9], [11], [13], [25]. The connection with the (deterministic
or random) abstract economy models has been intensively approached in [26].
Variants of abstract economy model (or generalized game) have been defined
by Shafer and Sonnenschein [20] or Yannelis and Prahbakar [24]. In [8] Kim
and Lee defined the free abstract economy and they also proved existence
theorems of best proximity pairs and equilibrium pairs which generalize the
previous results due to Srinivasan and Veeramani [21],[22], Sehgal and Singh
[19] or Reich [18]. Kim [7] also obtained generalizations of the theorems
from [8].

In this paper, we introduce new systems of generalized quasi-variational
inequalities and we prove the existence of the solutions by using results of
pair equilibrium existence for free abstract economies. We also propose a
fuzzy approach of this topic and we obtain results concerning the existence
of solutions for the random systems of generalized quasi-variational
inequalities.

The paper is organized in the following way: Section 2 contains
preliminaries and notation. New types of systems of generalized
quasi-variational inequalities are introduced in Section 3. The main results
are stated in Section 4. In Section 5 a fuzzy approach of the topic is
proposed. Finally, the systems of random quasi-variational inequalities are
studied in Section 6. A random fixed point theorem is obtained as a
corollary.

\section{Preliminaries and notation}

\subsection{Definitions and notation}

Throughout this paper, we shall use the following notation and definitions:

Let $A$ be a subset of a topological space $X$. 2$^{A}$ denotes the family
of all subsets of $A$. cl$A$ denotes the closure of $A$ in $X$. If $A$ is a
subset of a vector space, co$A$ denotes the convex hull of $A$. If $F$, $T:$ 
$A\rightarrow 2^{X}$ are correspondences, then co$T$, cl$T$, $T\cap F$ $:$ $%
A\rightarrow 2^{X}$ are correspondences defined by $($co$T)(x)=$co$T(x)$, $($%
cl$T)(x)=$cl$T(x)$ and $(T\cap F)(x)=T(x)\cap F(x)$ for each $x\in A$,
respectively. The graph of $T:X\rightarrow 2^{Y}$ is the set Gr$%
(T)=\{(x,y)\in X\times Y\mid y\in T(x)\}.\medskip $

\begin{notation}
Let $X$ and $Y$ be two nonempty sets. We denote by $\mathcal{F}(Y)$ the
collection of fuzzy sets on $Y$. A mapping from $X$ into $\mathcal{F}(Y)$ is
called a fuzzy mapping. If $F:X\rightarrow \mathcal{F}(Y)$ is a fuzzy
mapping, then for each $x\in X,$ $F(x)$ (denoted by $F_{x}$ in this sequel)
is a fuzzy set in $\mathcal{F}(Y)$ and $F_{x}(y)$ is the degree of
membership of point $y$ in $F_{x}.$
\end{notation}

Let E and F be two Hausdorff topological vector spaces and $X\subset E$, $%
Y\subset F$ be two nonempty convex subsets. A fuzzy mapping $F:X\rightarrow 
\mathcal{F}(Y)$ is called convex, if for each $x\in X,$ the fuzzy set $F_{x}$
on $Y$ is a fuzzy convex set, i.e., for any $y_{1},y_{2}\in Y,$ $t\in
\lbrack 0,1],$ $F_{x}(ty_{1}+(1-t)y_{2})\geq \min
\{F_{x}(y_{1}),F_{x}(y_{2})\}.$

In the sequel, we denote by

$(A)_{a}=\{y\in Y:A(y)\geq a\},$ $a\in \lbrack 0,1]$ the $a$-cut set of $%
A\in \mathcal{F}(Y).\medskip $

\subsection{Continuity of correspondences}

In this subsection we remind several notions concerning the continuity of
the correspondences defined on topological spaces. The motivation for our
presentation is the fact that the correspondences are the main mathematical
objects which define the notions from the main results concerning the free
abstract economies and variational inequalities.

\begin{definition}
Let $X$, $Y$ be topological spaces and $T:X\rightarrow 2^{Y}$ be a
correspondence
\end{definition}

\begin{enumerate}
\item $T$ is said to be \textit{upper semicontinuous} if for each $x\in X$
and each open set $V$ in $Y$ with $T(x)\subset V$, there exists an open
neighborhood $U$ of $x$ in $X$ such that $T(x)\subset V$ for each $y\in U$.

\item $T$ is said to be \textit{lower semicontinuous} (shortly l.s.c) if for
each x$\in X$ and each open set $V$ in $Y$ with $T(x)\cap V\neq \emptyset $,
there exists an open neighborhood $U$ of $x$ in $X$ such that $T(y)\cap
V\neq \emptyset $ for each $y\in U$.

\item $T$ is said to have \textit{open lower sections} if $T^{-1}(y):=\{x\in
X:y\in T(x)\}$ is open in $X$ for each $y\in Y.$

\item The correspondence $\overline{T}$ is defined by $\overline{T}%
(x)=\{y\in Y:(x,y)\in $cl$_{X\times Y}$Gr$(T)\}$ (the set cl$_{X\times Y}$Gr$%
(T)$ is called the adherence of the graph of T)$.$ It is easy to see that cl$%
T(x)\subset \overline{T}(x)$ for each $x\in X.$\medskip
\end{enumerate}

The following lemma helps us to construct new correspondences having
different properties of continuity.

\begin{lemma}
(see [26])\textit{Let }$X$\textit{\ and }$Y$\textit{\ be two topological
spaces. }
\end{lemma}

\textit{Suppose }$T_{1}:X\rightarrow 2^{Y}$\textit{\ , }$T_{2}:X\rightarrow
2^{Y}$\textit{\ are upper semicontinuous(resp. lower semicontinuous) such
that }$T_{2}(x)\subset T_{1}(x)$\textit{\ for all }$x\in A,$ where\textit{\ }%
$A$\textit{\ is an open (resp. closed) subset of }$X.$\textit{\ Then the
correspondence }$T:X\rightarrow 2^{Y}$\textit{\ defined by}%
\begin{equation*}
T(x)=\left\{ 
\begin{array}{c}
T_{1}(x)\text{, \ \ \ \ \ \ \ if }x\notin A\text{, } \\ 
T_{2}(x)\text{, \ \ \ \ \ \ \ \ \ if }x\in A%
\end{array}%
\right.
\end{equation*}%
\textit{is also upper semicontinuous (resp. lower semicontinuous).}

\begin{definition}
Let $X$, $Y$ be topological spaces and $T:X\rightarrow 2^{Y}$ be a
correspondence. An element $x\in X$ is named \textit{maximal} \textit{element%
} for $T$ if $T(x)=\Phi .$
\end{definition}

\textit{Let }$I$\textit{\ be an index set.} \textit{For each }$i\in I,$%
\textit{\ let }$X_{i}$\textit{\ be a nonempty subset of a topological space }%
$E_{i}$\textit{\ and }$T_{i}:X:=\underset{i\in I}{\prod }X_{i}\rightarrow
2^{Y_{i}}$\textit{\ a correspondence. Then a point }$x\in X$\textit{\ is
called a maximal element for the family of correspondences }$\{T_{i}\}_{i\in
I}$\textit{\ if }$T_{i}(x)=\emptyset $\textit{\ for all }$i\in I.$

\textit{The family }$\left( X_{i},T_{i}\right) _{i\in I}$\textit{\ is called
a qualitative game.}

The abstract economies are extensions of the qualitative games. We present
the definition below.

\begin{definition}
A \emph{generalized abstract economy} (or \emph{generalized game}) $\Gamma $
is defined as a family $(X_{i},A_{i},P_{i})_{i\in I}$ where $%
A_{i}:X\rightarrow 2^{X_{i}}$ is a constraint correspondence and $%
P_{i}:X\rightarrow 2^{X_{i}}$ is a preference correspondence$.$
\end{definition}

\begin{definition}
An \emph{equilibrium }for $\Gamma $ is a point $x^{\ast }\in X$ such that
for each $i\in I,$ $x_{i}^{\ast }\in $ $A_{i}(x^{\ast })$ and $P_{i}(x^{\ast
})\cap A_{i}(x^{\ast })=\emptyset .$
\end{definition}

Theorem 1 is a maximal element theorem for upper semicontinuous
correspondences. We will use it in order to prove the theorems concerning
the existence of equilibrium pairs for the abstract economies.

\begin{theorem}
(X. P. Ding, [5]) \textit{Let }$X$\textit{\ be a nonempty subset of a
locally convex Hausdorff topological vector space\ and }$D$\textit{\ a
nonempty, compact subset of }$X.$\textit{\ Let }$T$\textit{\ }$:X\rightarrow
2^{D\text{ }}$\textit{\ be an upper semicontinuous correspondence such that }%
for each $x\in X$, $x\notin $clco$T(x)$.\textit{\ Then there exists }$%
x^{\ast }\in $\textit{co}$D$\textit{\ such that }$T\left( x^{\ast }\right)
=\emptyset .$
\end{theorem}

Theorem 2 is an existence theorem for maximal elements that is Theorem 7 in
[23].

\begin{theorem}
(see [23])\textit{\ Let\ }$\Gamma =\left( X_{i},T_{i}\right) _{i\in I}$%
\textit{\ be a qualitative game where }$I$\textit{\ is an index set such
that for each }$i\in I,$ \textit{the following conditions hold:}
\end{theorem}

\textit{1) }$X_{i}$\textit{\ is a nonempty convex compact metrizable subset
of a Hausdorff locally convex topological vector space }$E$\textit{\ and }$%
X:=\underset{i\in I}{\prod }X_{i}$\textit{,}

\textit{2) }$T_{i}:X\rightarrow 2^{X_{i}}$\textit{\ is lower semi-continuous;%
}

4) \textit{for each} $x\in X,$ $x_{i}\notin $clco$T_{i}(x);$

\textit{Then there exists a point }$x^{\ast }\in X$\textit{\ such that }$%
T_{i}(x^{\ast })=\emptyset $\textit{\ for all }$i\in I,$ \textit{i.e. }$%
x^{\ast }$\textit{\ is a maximal element of }$\Gamma .$\textit{\bigskip }

The maximal element theorems can be consequences of the equilibrium
theorems, as it can be seen in the following situation.

We proved in [16] the next equilibrium theorem.

\begin{theorem}
\textit{Let }$\Gamma =(X_{i},A_{i},P_{i},B_{i})_{i\in I}$\textit{\ \ be an
abstract economy, where }$I$\textit{\ is a (possibly uncountable) set of
agents such that for each }$i\in I:$
\end{theorem}

(1)\textit{\ }$X_{i}$\textit{\ is a non-empty convex set in a Hausdorff
locally convex space }$E_{i}$\textit{, }$X:=\underset{i\in I}{\prod }X_{i}$%
\textit{\ \ is paracompact and }$D_{i}$ is \textit{a non-empty, convex,
compact subset of }$X_{i}$\textit{;}

(2)\textit{\ }$B_{i}$\textit{\ is lower semicontinuous with non-empty convex
values and } \textit{for each }$x\in X$,\textit{\ }$A_{i}(x)\neq \emptyset $%
, $A_{i}(x)\subset B_{i}(x)\ $\textit{and }cl$B_{i}(x)\cap D_{i}\neq
\emptyset $\textit{;}

(3) \textit{\ the correspondence }$A_{i}$\textit{\ }$\cap P_{i}:X\rightarrow
2^{D_{i}}$ \textit{is lower semi-continuous;}

(4) \textit{for each} $x\in X,$ $x_{i}\notin \overline{(\text{co}A_{i}\cap 
\text{co}P_{i})}(x).$

\textit{Then there exists an equilibrium point }$x^{\ast }\in D$ \textit{\
for }$\Gamma $\textit{,}$\ i.e.$\textit{, for each }$i\in I$\textit{, }$%
x_{i}^{\ast }\in \overline{B}_{i}(x^{\ast })$\textit{\ and }$A_{i}(x^{\ast
})\cap P_{i}(x^{\ast })=\emptyset $\textit{.\medskip }

Theorem 4 is an existence theorem for maximal elements that is a consequence
of Theorem 3.

\begin{theorem}
\textit{Let\ }$\Gamma =\left( X_{i},T_{i}\right) _{i\in I}$\textit{\ be a
qualitative game where }$I$\textit{\ is an index set such that for each }$%
i\in I,$ \textit{the following conditions hold:}
\end{theorem}

\textit{1) }$X_{i}$\textit{\ is a nonempty convex compact subset of a
Hausdorff locally convex topological vector space }$E$\textit{\ and }$X:=%
\underset{i\in I}{\prod }X_{i}$\textit{,}

\textit{2) }$T_{i}:X\rightarrow 2^{X_{i}}$\textit{\ is lower semi-continuous;%
}

3) \textit{for each} $x\in X,$ $x_{i}\notin $co$\overline{T}_{i}(x)$

\textit{Then there exists a point }$x^{\ast }\in X$\textit{\ such that }$%
T_{i}(x^{\ast })=\emptyset $\textit{\ for all }$i\in I,$ \textit{i.e. }$%
x^{\ast }$\textit{\ is a maximal element of }$\Gamma .$

\subsection{Best proximity pairs of correspondences}

This subsection traits the problem of existence of best proximity pairs for
correspondences defined on normed linear spaces. The best proximity pair
theorems analyze the conditions under which the problem of minimizing the
real-valued function $x\rightarrow d$($x,T\left( x\right) )$ has a solution.

Firstly, we present the following notation.

\begin{notation}
Let $X$ and $Y$ be any two subsets of a normed space $E$ with a norm $%
\parallel \cdot \parallel $, and the metric $d(x,y)$ is induced by the norm.
Then, we now recall the following notation:
\end{notation}

Prox$(X,Y):=\{(x,y)\in X\times Y:d(x,y)=d(X,Y)=$inf$\{d(x,y):x\in X,y\in
Y\}\};$

$X_{0}:=\{x\in X:d(x,y)=d(X,Y)$ for some $y\in Y\};$

$Y_{0}:=\{y\in Y:d(x,y)=d(X,Y)$ for some $x\in X\}.$

If $X$ and $Y$ are non-empty compact and convex subsets of a normed linear
space, then it is easy to see that $X_{0}$ and $Y_{0}$ are both non-empty
compact and convex.

Let $I$ be a finite (or an infinite) index set. For each $i\in I,$ let $X$
and $Y_{i}$ be a nonempty subsets of a normed space $E$ with a norm $%
\parallel \cdot \parallel $, and the metric $d(x,y)$ is induced by the norm.
Then, we can use the following notation: for each $i\in I,$

$X^{0}:=\{x\in X:$ for each $i\in I,$ $\exists $ $y_{i}\in Y_{i}$ such that $%
d(x,y_{i})=d(X,Y_{i})=$inf$\{d(x,y):x\in X,y\in Y_{i}\}\};$

$Y_{i}^{0}:=\{y\in Y_{i}:$ there exists $x\in X$ such that $%
d(x,y)=d(X,Y_{i})\}.$

When $\mid I\mid =1,$ it is easy to see that $X_{0}=X^{0}$ and $%
Y_{0}=Y_{i}^{0}.$

To prove our equilibrium theorems we need the following results.

\begin{definition}
(see [8])Let X and Y be two non-empty subsets of a normed linear space E,
and let $T:X\rightarrow 2^{Y}$ be a correspondence. Then the pair ($\mathit{x%
}^{\ast },T\left( \mathit{x}^{\ast }\right) $ is called the \textit{best
proximity pair} [8] for $T$ if $d$($\mathit{x}^{\ast },T\left( \mathit{x}%
^{\ast }\right) )=d$($\mathit{x}^{\ast },\mathit{y}^{\ast })=d(X,Y)$ for
some $\mathit{y}^{\ast }\in T(\mathit{x}^{\ast }).$
\end{definition}

The best proximity pair theorem seeks an appropriate solution which is
optimal.

W. K. Kim and K. H. Lee [8] gave the following result of existence of best
proximity pairs. This theorem is widely used in order to prove the existence
of the equilibrium pairs for the free abstract economies, and then, the
existence of the solutions for systems of variational quasi-inequalities.

\begin{theorem}
(see [8])For each $i\in I=\{1,...n\},$ let $X$ and $Y$ be non-empty compact
and convex subsets of a normed linear space $E$, and let $T_{i}:X\rightarrow
2^{Y_{i}}$ be an upper semicontinuous correspondence in $X^{0}$ such that $%
T_{i}\left( x\right) $ is nonempty closed and convex subset of $Y_{i}$ for
each $x\in X.$ Assume that $T_{i}\left( x\right) \cap Y_{i}^{0}\neq
\emptyset $ for each $x\in X^{0}.$
\end{theorem}

Then there exists a system of best proximity pairs $\left\{ \mathit{x}^{\ast
}\right\} \times T_{i}\left( \mathit{x}^{\ast }\right) \subseteq X\times
Y_{i},i.e.,$ for each $i\in I,$ $d(\mathit{x}^{\ast },T_{i}\left( \mathit{x}%
^{\ast }\right) )=d(X,Y_{i}).$

\begin{definition}
(see [8]). The set $\mathcal{A}_{x}=\{y\in Y:y\in A(x)$ and $d(x,y)=d(X,Y)\}$
is named \textit{the best proximity set} of the correspondence $%
A:X\rightarrow 2^{Y}$ at $x$.
\end{definition}

In general, $\mathcal{A}_{x}$ might be an empty set. If $(x,A(x))$ is a
proximity pair for $A$ and $A(x)$ is compact, then $\mathcal{A}_{x}$ must be
non-empty.

\section{New types of systems of generalized quasi-variational inequalities}

In this section, we introduce some new types of system of generalized
quasi-variational inequalities.

Let $I$ be a finite set. For each $i\in I,$ let $X$ and $Y_{i}$ be non-empty
compact and convex subsets of a normed linear space $E$, $E^{\prime }$ be
the dual space of $E$, let $A_{i}:X\rightarrow 2^{Y_{i}},$ $%
B_{i}:Y_{i}\rightarrow 2^{E^{\prime }}$ be correspondences and $\psi
_{i}:Y\times Y_{i}\rightarrow \mathbb{R}\cup \{-\infty ,+\infty \}$.

We associate with $A_{i}$ and $G_{i}$ the next generalized quasi-variational
problems:

(1) Find the pair $(x^{\ast },y^{\ast })=(x^{\ast },(y_{i}^{\ast })_{i\in
I})\in X\times \tprod_{i\in I}Y_{i}$ such that:\medskip

i)$y_{i}^{\ast }\in A_{i}(x);$

ii)$d(x^{\ast },y_{i}^{\ast })=d(X,Y_{i});$

iii)$\sup_{z_{i}\in A_{i}(x^{\ast })}\psi _{i}(y^{\ast },z_{i})\leq 0.$%
\medskip

If we consider a particular function $\psi _{i},$ we have the following
problem:

(2) Find the pair $(x^{\ast },y^{\ast })=(x^{\ast },(y_{i}^{\ast })_{i\in
I})\in X\times \tprod_{i\in I}Y_{i}$ such that:\medskip

i) $y_{i}^{\ast }\in A_{i}(x^{\ast });$

ii) $d(x^{\ast },y_{i}^{\ast })=d(X,Y_{i});$

iii) $\sup_{z_{i}\in A_{i}(x^{\ast })}\sup_{v\in G_{i}(y_{i}^{\ast })}\func{%
Re}\langle v,y_{i}^{\ast }-z_{i}\rangle \leq 0;$

where the real part of pairing between $E^{\prime }$ and $E$ is denoted by $%
\func{Re}\langle v,x\rangle $ for each $v\in E$ and $x\in E.$

These systems of generalized quasi-variational inequalities generalize that
ones studied, for instance, by Yuan in [26].\medskip 

We will also work in the fuzziness framework.

Let $I$ be a nonempty set. For each $i\in I$, let $X_{i}$ and $Y_{i}$ be
non-empty subsets of a normed linear space $E$. Define $X:=\underset{i\in I}{%
\prod }X_{i}$; let $A_{i}:X\rightarrow \mathcal{F}(Y_{i})$ be a fuzzy
correspondence, $P_{i}:Y:=\underset{i\in I}{\prod }Y_{i}\rightarrow \mathcal{%
F}(Y_{i})$ a fuzzy mapping and $a_{i}:X\rightarrow (0,1]$ a fuzzy function.
We consider that $X_{i}$ and $Y_{i}$ are non-empty subsets of a normed
linear space $E$.

Let $G_{i}:Y_{i}\rightarrow 2^{E_{i}^{\prime }}$ be a correspondence and $%
\psi _{i}:Y\times Y_{i}\rightarrow \mathbb{R}\cup \{-\infty ,+\infty \}$.

Whenever $X_{i}=X$ for each $i\in I,$ for the simplicity, we may assume $%
A_{i}:X\rightarrow $ $\mathcal{F}(Y_{i})$ instead of $\ A_{i}:\underset{i\in
I}{\prod }X_{i}\rightarrow \mathcal{F}(Y_{i}).$

Now, we are introducing the next type of system of quasi-variational
inequalities:

(3) Find the pair $(x^{\ast },y^{\ast })\in X\times Y$ such that for every%
\textit{\ }$i\in I$\textit{:}

i) $y_{i}^{\ast }\in (A_{i}(x^{\ast }))_{a_{i}(x^{\ast })};$

ii) $d(x^{\ast },y_{i}^{\ast })=d(X,Y_{i});$

iii) $\sup_{z_{i}\in (A_{i}(x^{\ast }))_{a_{i}(x^{\ast })}}\psi _{i}(y^{\ast
},z_{i})\leq 0,$

where $(A_{i_{x^{\ast }}})_{a_{i}(x^{\ast })}=\{z\in Y_{i}:A_{i_{x^{\ast
}}}(z)\geq a_{i}(x^{\ast })\}.$\medskip

If $A_{i}:X\rightarrow 2^{Y_{i}}$ is a classical correspondence, then we get
the system of quasi-variational inequalities defined at (1).

We will also work with the next fuzzy model:

(4) Find the pair $(x^{\ast },y^{\ast })\in X\times Y$ such that for every%
\textit{\ }$i\in I$\textit{:}

i) $y_{i}^{\ast }\in (A_{i}(x^{\ast }))_{a_{i}(x^{\ast })};$

ii) $d(x^{\ast },y_{i}^{\ast })=d(X,Y_{i});$

iii) $\sup_{z_{i}\in (A_{i}(x^{\ast }))_{a_{i}(x^{\ast })}}\sup_{v\in
G_{i}(y_{i}^{\ast })}\func{Re}\langle v,y_{i}^{\ast }-z_{i}\rangle \leq 0$%
.\medskip

Noor and Elsanousi [15] introduced the notion of a random variational
inequality.

We propose the next random system of quasi-variational inequalities which
generalizes the random one studied, for instance, by Yuan in [26].

Let $(\Omega ,\mathcal{F})$ be a measurable space. For each $i\in I$, let $%
X_{i}$ and $Y_{i}$ be non-empty subsets of a normed linear space $E$. Define 
$X:=\underset{i\in I}{\prod }X_{i}$; let $A_{i}:\Omega \times X\rightarrow
2^{Y_{i}}.$

(5) Find the measurable functions $\varphi ^{1}:\Omega \rightarrow X$ and $%
\varphi ^{2}:\Omega \rightarrow Y$ such that for each $i\in I$ and for all $%
\omega \in \Omega :$

i) $\pi _{i}(\varphi ^{2}(\omega ))\in A_{i}(\omega ,\varphi ^{1}(\omega ));$

ii) $d(\varphi ^{1}(\omega ),\pi _{i}(\varphi ^{2}(\omega )))=d(X,Y_{i});$

iii) $\sup_{z_{i}\in A_{i}(\varphi ^{1}(\omega ))}\psi _{i}(\varphi
^{2}(\omega ),z_{i})\leq 0.$ \textit{\medskip }

\section{Main results}

The aim of this section is to state theorems concerning the pair equilibrium
existence for free abstract economies with Q-majorized preference
correspondences and to apply them in order to prove the existence of the
solutions for the systems of quasi-variational inequalities introduced in
Section 3.

First, we present the model of a free abstract economy introduced by Kim and
Lee [8].

Let $I$ be a nonempty set (the set of agents). For each $i\in I$, let $X_{i}$
be a non-empty set of manufacturing commodities, and $Y_{i}$ be a non-empty
set of selling commodities. Define $X:=\underset{i\in I}{\prod }X_{i}$; let $%
A_{i}:X\rightarrow 2^{Y_{i}}$ be the constraint correspondence and $P_{i}:Y:=%
\underset{i\in I}{\prod }Y_{i}\rightarrow 2^{Y_{i}}$ the preference
correspondence. We consider that $X_{i}$ and $Y_{i}$ are non-empty subsets
of a normed linear space $E$.

\begin{definition}
\textit{\ }A \textit{free} \textit{abstract economy} is the family $\Gamma
=(X_{i},Y_{i},A_{i},P_{i})_{i\in I}$.\medskip
\end{definition}

\begin{definition}
\textit{\ }An \textit{equilibrium} pair for $\Gamma $ is defined as a pair
of points $(x^{\ast },y^{\ast })\in X\times Y$ such that for each $i\in I$, $%
y_{i}^{\ast }\in B_{i}(x^{\ast })$ with $d(x_{i}^{\ast },y_{i}^{\ast
})=d(X_{i},Y_{i})$ and $A_{i}(x^{\ast })\cap P_{i}(y^{\ast })=\emptyset $%
.\medskip
\end{definition}

Whenever $X_{i}=X$ for each $i\in I,$ for the simplicity, we may assume $%
A_{i}:X\rightarrow 2^{Y_{i}}$ instead of $\ A_{i}:\underset{i\in I}{\prod }%
X_{i}\rightarrow 2^{Y_{i}}$ for the free abstract economy $\Gamma
=(X,Y_{i},A_{i},P_{i})_{i\in I}$ and equilibrium pair. In particular, when $%
I=\{1,2...n\},$ we may call $\Gamma $ a free n-person game.

The economic interpretation of an \textit{equilibrium} pair for $\Gamma $ is
based on the requirement that for each $i\in I,$ minimize the travelling
cost $d(x_{i},y_{i})$, and also, maximize the preference $P_{i}(y)$ on the
constraint set $A_{i}(x)$. Therefore, it is contemplated to find a pair of
points $(x^{\ast },y^{\ast })\in X\times Y$ such that for each $i\in I$, $%
y_{i}^{\ast }\in B_{i}(x^{\ast })$, $A_{i}(x^{\ast })\cap P_{i}(y^{\ast
})=\emptyset $ and $\parallel x_{i}^{\ast }-y_{i}^{\ast }\parallel
=d(X_{i},Y_{i}),$ where $d(X_{i},Y_{i})=\inf \left\{ \parallel x_{i}^{\ast
}-y_{i}^{\ast }\parallel \mid x_{i}\in X_{i},\text{ }y_{i}\in Y_{i}\right\} $%
.

When in addition $X_{i}=Y_{i}$ for each $i\in I,$ then the previous
definitions can be reduced to the standard definitions of equilibrium theory
in economics due to Shafer and Sonnenshein [20] or Yannelis and Prabhakar
[24].\medskip 

In order to prove our main theorems, we need the following results
concerning the $Q-$ majorized correspondences, which generalize the lower
semicontinuous ones.

\begin{definition}
(see [12])Let $X$ be a topological space and $Y$ be a non-empty subset of a
vector space $E$, $\theta :X\rightarrow E$ a function and $\,T:X\rightarrow
2^{Y}$ a correspondence.{}
\end{definition}

1) $T$ \textit{is of class }$Q_{\theta }$ (or $Q$) if:

\qquad a) for each $x\in X$, $\theta (x)\notin $cl$T(x)$ and

\qquad b) $T$ is lower semicontinuous with open and convex values in $Y$;

2) A correspondence $T_{x}:X\rightarrow 2^{Y}$ is a $Q_{\theta }$\textit{%
-majorant of }$T$ at $x,$ \thinspace if there exists an open neighborhood $%
N(x)$ of $x$ such that:

\qquad a) for each $z\in N(x)$, $T(z)\subset T_{x}(z)$ and $\theta (z)\notin 
$cl$T_{x}(z);$

\qquad b) $T_{x}$ is lower semicontinuous with open convex values;

3) $T$\textit{\ is }$Q_{\theta }$\textit{-majorized} if for each $x\in X$
with $T(x)\neq \emptyset ,$ there exists a $Q_{\theta }$-majorant $T_{x}$ of 
$T$ at $x$.\medskip

The next property is remarkable for the $Q$-majorized correspondences. It
says that a $Q$-majorized correspondence defined on a regular paracompact
topological vector space is a selector of a correspondence of class $%
Q_{\theta }$ defined on the whole space.\medskip

\begin{theorem}
(see [12]) \textit{Let }$X$\textit{\ be a regular paracompact topological
vector space and }$Y$\textit{\ be a non-empty subset of \ a vector space }$E$%
\textit{. Let }$\theta :X\rightarrow E$\textit{\ and }$T:X\rightarrow
2^{Y}\smallsetminus \left\{ \emptyset \right\} $\textit{\ be a }$Q_{\theta }-
$\textit{majorized} \textit{correspondence. Then, there exists a
correspondence }$S:X\rightarrow 2^{Y}$\textit{\ of class }$Q_{\theta }$%
\textit{\ such that }$T(x)\subset S(x)$\textit{\ for each }$x\in X.\medskip $
\end{theorem}

Theorem 7 is an existence theorem of pair equilibrium for a free n person
game with upper semi-continuous constraint correspondences and $Q_{\theta }$%
-majorized preference correspondences.

\begin{theorem}
\textit{Let }$\Gamma =(X,Y_{i},A_{i},P_{i})_{i\in I}$\textit{\ be a free
n-person game such that for each }$i\in I=\{1,2...n\}$:
\end{theorem}

(1)\textit{\ }$X$\textit{\ and }$\mathit{Y}_{i}$\textit{\ are non-empty
compact and convex subsets of normed linear space }$\mathit{E}$\textit{;}

(2)\textit{\ }$A_{i}:X\rightarrow 2^{Y_{i}}$\textit{\ is such that }$A_{i}$%
\textit{\ is upper semicontinuous in }$X^{0},$\textit{\ }$A_{i}(x)$\textit{\
is a closed convex subset of }$Y_{i}$\textit{, }$A_{i}(x)\cap Y_{i}^{0}\neq
\emptyset $\textit{\ for each }$x\in X^{0}$\textit{;}

(3)\textit{\ }$P_{i}:Y:\underset{i\in I}{=\prod }Y_{i}\rightarrow 2^{Y_{i}}$ 
\textit{is }$Q_{\pi _{i}}-$\textit{majorized}$;$

(4) $P_{i}(y)$ \textit{is nonempty for each} $y\in Y$;

\textit{Then there exists an equilibrium pair} \textit{of points} $(x^{\ast
},y^{\ast })\in X\times Y$ \textit{such that for each }$i\in I$\textit{, }$%
y_{i}^{\ast }\in P_{i}(x^{\ast })$\textit{\ with} $d(x^{\ast },y_{i}^{\ast
})=d(X,Y_{i})$ \textit{and} $A_{i}(x^{\ast })\cap P_{i}(y^{\ast })=\emptyset 
$.\textit{\medskip }

\textit{Proof.} Since $A_{i}$ satisfies the whole assumption of Theorem 5
for each $i\in I$, there exists a point $x^{\ast }\in X$ satisfying the
system of best proximity pairs, i.e., for each $i\in I,$ $\left\{ x^{\ast
}\right\} \times A_{i}(x^{\ast })\subseteq X\times Y_{i}$ such that $%
d(x^{\ast },A_{i}(x^{\ast }))=d(X,Y_{i}).$ Let $\mathcal{A}_{i}:=\left\{
y\in A_{i}(x^{\ast })\text{ / }d(x^{\ast },y)=d(X,Y_{i})\right\} $ the
non-empty best proximity set of the correspondence $A_{i}.$ The set $%
\mathcal{A}_{i}$ is nonempty, closed and convex.

Since $P_{i}$ is $Q_{\pi _{i}}$\textit{-}majorized for each $i\in I$, by
Theorem 6, there exists a correspondence $\varphi _{i}:Y\rightarrow
2^{Y_{i}} $ of class $Q_{\pi _{i}}$ such that $P_{i}(y)\subset \varphi
_{i}(y)$ for each $y\in Y$ and for each $i\in I.$ Then, $\varphi _{i}$ is
lower semicontinuous with open, convex values and $\pi _{i}(y)\notin $cl$%
\varphi _{i}(y).$

For each $i\in I$ define the correspondence

$\Phi _{i}:Y\rightarrow 2^{Y_{i}}$ by

$\Phi _{i}(y):=\left\{ 
\begin{array}{c}
\varphi _{i}(y),\text{ \ \ \ \ \ if }y_{i}\notin \mathcal{A}_{i}, \\ 
A_{i}(x^{\ast })\cap \varphi _{i}(y),\text{ \ if }y_{i}\in \mathcal{A}_{i}.%
\end{array}%
\right. $

By Lemma 1, $\Phi _{i}$ is lower semicontinuous, has convex values, and $\pi
_{i}(y)\notin $cl$\Phi _{i}(y).$ By applying Theorem $4$ to $(Y_{i},\Phi
_{i})_{i\in I}$, there exists a maximal element $y^{\ast }\in Y$ such that $%
\Phi _{i}\left( y^{\ast }\right) =\emptyset $ for each $i\in I.$ Since $%
P_{i}(y)\neq \emptyset $ for each $y\in Y,$ $\Phi _{i}(y)$ is a non-empty
subset of $Y$ for each $y\in Y$ with $y_{i}\notin \mathcal{A}_{i}$. It
follows that $y_{i}^{\ast }\in \mathcal{A}_{i}$ and $A_{i}(x^{\ast })\cap
\varphi _{i}(y^{\ast })=\emptyset $. We have that $A_{i}(x^{\ast })\cap
P_{i}(y^{\ast })=\emptyset $ because $P_{i}(y^{\ast })\subset \varphi
_{i}(y^{\ast }).$ Hence, $y_{i}^{\ast }\in \mathcal{A}_{i}$ such that $%
d(x^{\ast },y_{i}^{\ast })=d(X,Y_{i})$ for each $i\in I$ and then $(x^{\ast
},y^{\ast })$ is an equilibrium pair for $\Gamma .$ $\square \medskip $

By using the theorem above, we obtain the following result concerning the
existence of the solutions for the systems of quasi-variational inequalities
of type (1), where the correspondences $A_{i}$ are upper semi-continuous.

\begin{theorem}
\textit{Let }$X$\textit{\ and }$Y_{i}$\textit{\ be non-empty compact and
convex subsets of a normed linear space }$E$\textit{\ for each }$i\in
I=\{1,2,...,n\}.$ Suppose that for each $i\in I,$ the following conditions
are fulfilled:
\end{theorem}

\textit{(1) }$A_{i}:X\rightarrow 2^{Y_{i}}$\textit{\ is an upper
semicontinuous correspondence such that each }$A_{i}(x)$\textit{\ is a
non-empty closed and convex subset of }$Y_{i}$\textit{\ and }$A_{i}(x)\cap
Y_{i}^{0}\neq \emptyset $ for each $x\in X^{0}\mathit{;}$

\textit{(2) The function\ }$\psi _{i}:Y\times Y_{i}\rightarrow R\cup
\{-\infty ,+\infty \}$\textit{\ is such that }

\ \ \ (\textit{a) for each fixed }$y_{i}\in Y_{i},$ $\psi _{i}(\cdot ,y_{i})$
\textit{is lower semicontinuous such that} $\{z\in Y_{i}:\psi _{i}(y,z)>0\}$ 
\textit{is non-empty for each }$y\in Y$\textit{;}

\ \ \ (\textit{b) }$y_{i}\notin \{z\in Y_{i}:\psi _{i}(y,z)>0\}$\textit{\
for each fixed }$y\in Y;$

\ \ \ (\textit{c) for each} $y\in Y,$ $\psi _{i}(y,\cdot )$ \textit{is
concave.}

\textit{Then, there exists }$(x^{\ast },y^{\ast })\in X\times Y$\textit{\
such that for every }$i\in I$\textit{:}

i) $y_{i}^{\ast }\in A_{i}(x^{\ast });$

ii) $d(x^{\ast },y_{i}^{\ast })=d(X,Y_{i});$

iii) $\sup_{z_{i}^{\ast }\in A_{i}(x^{\ast })}\psi _{i}(y^{\ast
},z_{i}^{\ast })\leq 0.$

\textit{Proof. }For every $i\in I,$ let $P_{i}:Y\rightarrow Y_{i}$ be
defined by $P_{i}(y)=\{z\in Y_{i}:\psi _{i}(y,z)>0\}$ for each $y\in Y.$

We shall show that the free abstract economy $G=(X,Y_{i},A_{i},P_{i},)_{i\in
I}$ satisfies all hypotheses of Theorem 7.

According to 2 a), we have that\textit{\ }$P_{i}$\textit{\ }has open lower
sections (and then lower semicontinuous) with nonempty compact values and
according to 2 b), $y_{i}\not\in P_{i}(y)$ for each $y\in Y.$ Assumption 2
c) implies that $P_{i}$ has convex values. Hence, $P_{i}$ is $Q-$majorized.

Thus the free abstract economy $G=(X,Y_{i},A_{i},P_{i})_{i\in I}$ satisfies
all hypotheses of Theorem 7. Therefore, there exists $(x^{\ast },y^{\ast
})\in X\times Y$ such that for every $i\in I:$ $y_{i}^{\ast }\in
A_{i}(x^{\ast })$; $A_{i}(x^{\ast })\cap P_{i}(y^{\ast })=\phi $ and $%
d(x^{\ast },y_{i}^{\ast })=d(X,Y_{i}),$ that is, there exists $(x^{\ast
},y^{\ast })\in X\times Y$ such that for every $i\in I:$

i) $y_{i}^{\ast }\in A_{i}(x^{\ast })$;

ii) sup$_{z\in A_{i}(x^{\ast })}\psi _{i}(y^{\ast },z)\leq 0$;

iii) $d(x^{\ast },y_{i}^{\ast })=d(X,Y_{i}).$ $\square \medskip $

Now, we are studying the case when the correspondences $A_{i}$, $i\in I$ are
lower semicontinuous. Firstly, we will prove the existence theorem of the
equilibrium pairs for a free abstract economy and then, we will used it in
order to prove a theorem which states the existence of the solutions of a
system of quasi-variational inequalities of type (1) with lower
semicontinuous correspondences $A_{i},$ $i\in I.$

Let $I$ be an index set. Suppose that for each $i\in I,$ $X$, $Y_{i}$\ are
non-empty compact and convex subsets of normed linear space $E$ and $%
A_{i}:X\rightarrow 2^{Y_{i}}$\textit{\ }is such that $A_{i}$\ is lower
semicontinuous in $X^{0}$\ and $A_{i}(x)$\ is a closed convex subset of $%
Y_{i}.$ By Theorem 1.1 in Michael [14], for each $i\in I$ there exists an
upper semicontinuous correspondence $H_{i}:X\rightarrow 2^{Y_{i}}$ with
nonempty values such that $H_{i}(x)\subset A_{i}(x)$ for all $x\in X.$

\begin{definition}
We say that $A_{i}$ has the property * if, in addition, $H_{i}(x)\cap
Y_{i}^{0}\neq \emptyset $\textit{\ }for each\textit{\ }$x\in X^{0}.$
\end{definition}

Theorem 9 is an existence theorem of pair equilibrium for a free n person
game with lower semi-continuous constraint correspondences and $Q_{\theta }$%
-majorized preference correspondences.

\begin{theorem}
Let $I=\{1,2,...,n\}$ and \textit{let }$\Gamma =(X,Y_{i},A_{i}P_{i})_{i\in
I} $\textit{\ be a free n-person game such that for each }$i\in I$\textit{:}
\end{theorem}

\textit{(1)\ }$X$\textit{\ and }$Y_{i}$\textit{\ are non-empty compact and
convex subsets of normed linear space E;}

\textit{(2)\ }$A_{i}:X\rightarrow 2^{Y_{i}}$\textit{\ is such that }$A_{i}$%
\textit{\ is lower semicontinuous in }$X^{0},$\textit{\ it has the property
* and }$A_{i}(x)$\textit{\ is a non-empty closed convex subset of }$Y_{i};$

\textit{(3)\ }$P_{i}:Y:\underset{i\in I}{=\prod }Y_{i}\rightarrow 2^{Y_{i}}$%
\textit{\ is }$Q_{\pi _{i}}-$\textit{majorized}$;$

\textit{(4) }$P_{i}(y)$\textit{\ is nonempty for each }$y\in Y$\textit{;}

\textit{Then there exists an equilibrium pair of points }$(x^{\ast },y^{\ast
})\in X\times Y$\textit{\ such that for each }$i\in I$\textit{, }$%
y_{i}^{\ast }\in A_{i}(x^{\ast })$\textit{\ with }$d(x^{\ast },y_{i}^{\ast
})=d(X,Y_{i})$\textit{\ and }$A_{i}(x^{\ast })\cap P_{i}(y^{\ast
})=\emptyset $\textit{.\medskip }

\textit{Proof.} Since $A_{i}$ has the property *, for each $i\in I$ there
exists an upper semicontinuous correspondence $H_{i}:X\rightarrow 2^{Y_{i}}$
with nonempty values such that $H_{i}(x)\subset B_{i}(x)$ for all $x\in X$
and $H_{i}(x)\cap Y_{i}^{0}\neq \emptyset $\textit{\ }for each\textit{\ }$%
x\in X^{0}.$ Let $S_{i}(x)=$clco$H_{i}(x)\subset A_{i}(x).$ The
correspondence $S_{i}$ satisfies the hypotheses of Theorem 5, then we get a
best proximity pair $\left\{ x^{\ast }\right\} \times S_{i}(x^{\ast
})\subseteq X\times Y_{i}$ for $S_{i}$, i.e. $d(x^{\ast },S_{i}(x^{\ast
}))=d(X,Y_{i}).$ Let $\mathcal{S}_{i}:=\left\{ y_{i}\in S_{i}(x^{\ast })%
\text{ / }d(x^{\ast },y_{i})=d(X,Y_{i})\right\} $ the non-empty best
proximity set of the correspondence $S_{i}.$ The set $\mathcal{S}_{i}$ is
nonempty, closed and convex.

Since $P_{i}$ is $Q_{\theta }$\textit{-}majorized, by Theorem $6$, we have
that there exists a correspondence $\varphi _{i}:Y\rightarrow 2^{Y_{i}}$ of
class $Q_{\pi _{i}}$ such that $P_{i}(y)\subset \varphi _{i}(y)$ for each $%
y\in Y.$ Then, $\varphi _{i}$ is lower semicontinuous with open, convex
values and $\pi _{i}(y)\notin $cl$\varphi _{i}(y).$

For each $i\in I$ define the correspondence $\Phi _{i}:Y\rightarrow
2^{Y_{i}} $ by

$\Phi _{i}(y):=\left\{ 
\begin{array}{c}
\varphi _{i}(y),\text{ \ \ \ \ \ if }y_{i}\notin \mathcal{S}_{i}, \\ 
A_{i}(x^{\ast })\cap \varphi _{i}(y),\text{ \ if }y_{i}\in \mathcal{S}_{i}.%
\end{array}%
\right. $

By Lemma 1, $\Phi _{i}$ is lower semicontinuous, it also has convex values
and $\pi _{i}(y)\notin $cl$\Phi _{i}(y).$ By applying Theorem $4$ to $%
(Y_{i},\Phi _{i})_{i\in I}$, there exists a maximal element $y^{\ast }\in X$
such that $\Phi _{i}(y^{\ast })=\emptyset $ for all $i\in I.$ Since $%
P_{i}(y)\neq \emptyset $ for each $y\in Y,$ $\Phi _{i}(y)$ is a non-empty
subset of $Y$ for each $y\in Y$ with $y_{i}\notin \mathcal{S}_{i}$. It
follows that $y_{i}^{\ast }\in \mathcal{S}_{i}$ and $A_{i}(x^{\ast })\cap
\varphi _{i}(y^{\ast })=\emptyset $. We have that $A_{i}(x^{\ast })\cap
P_{i}(y^{\ast })=\emptyset $ because $P_{i}(y^{\ast })\subset \varphi
_{i}(y^{\ast }).$ Hence, $y_{i}^{\ast }\in S_{i}(x^{\ast })\subset
A_{i}(x^{\ast })$ such that $d(x^{\ast },y_{i}^{\ast })=d(X,Y_{i})$ for each 
$i\in I$ and then $(x^{\ast },y^{\ast })$ is an equilibrium pair for $\Gamma
.$ $\square \medskip $

We are obtaining the following theorem concerning the systems of generalized
quasi-variational inequalities with lower semicontinuos correspondences $%
A_{i},$ $i\in I$.

\begin{theorem}
\textit{Let }$X$\textit{\ and }$Y_{i}$\textit{\ be non-empty compact and
convex subsets of a normed linear space }$E$\textit{\ for each }$i\in
I=\{1,2,...,n\}.$ Assume that for each $i\in I,$ the following conditions
are fulfilled:
\end{theorem}

\textit{The correspondence} $A_{i}:X\rightarrow 2^{Y_{i}}$\textit{\ is such
that}

\textit{(1) }$A_{i}$\textit{\ is lower semicontinuous in }$X^{0},$\textit{\
it has the property * and }$A_{i}(x)$\textit{\ is a non-empty closed convex
subset of }$Y_{i};$

\textit{The function\ }$\psi _{i}:Y\times Y_{i}\rightarrow \mathbb{R}\cup
\{-\infty ,+\infty \}$\textit{\ is such that }

\ \ \ (\textit{2) for each fixed }$y_{i}\in Y_{i},$ $\psi _{i}(\cdot ,y_{i})$
\textit{is lower semicontinuous and} $\{z\in Y_{i}:\psi _{i}(y,z)>0\}$ 
\textit{is non-empty for each }$y\in Y$\textit{;}

\ \ \ (\textit{3) }$y_{i}\notin \{z\in Y_{i}:\psi _{i}(y,z)>0\}$\textit{\
for each fixed }$y\in Y;$

\ \ \ (\textit{4) for each} $y\in Y,$ $\psi _{i}(y,\cdot )$ \textit{is
concave.}

\textit{Then, there exists }$(x^{\ast },y^{\ast })\in X\times Y$\textit{\
such that for every }$i\in I$\textit{:}

i) $y_{i}^{\ast }\in A_{i}(x^{\ast });$

ii) $d(x^{\ast },y_{i}^{\ast })=d(X,Y_{i});$

iii) $\sup_{z_{i}\in A_{i}(x^{\ast })}\psi _{i}(y^{\ast },z_{i})\leq 0.$

\textit{Proof. }For every $i\in I,$ let $P_{i}:Y\rightarrow Y_{i}$ be
defined by $P_{i}(y)=\{z\in Y_{i}:\psi _{i}(y,z)>0\}$ for each $y\in Y.$

We shall show that the free abstract economy $G=(X,Y_{i},A_{i},P_{i},)_{i\in
I}$ satisfies all hypotheses of Theorem 9.

According to 2), we have that\textit{\ }$P_{i}$\textit{\ }is lower
semicontinuous with nonempty compact values and according to 3), $%
y_{i}\not\in P_{i}(y)$ for each $y\in Y.$ Assumption 4) implies that $P_{i}$
has convex values. Then, $P_{i}$ is $Q$-majorized.

Thus the free abstract economy $G=(X,Y_{i},A_{i},P_{i})_{i\in I}$ satisfies
all hypotheses of Theorem 9. Therefore, there exists $(x^{\ast },y^{\ast
})\in X\times Y$ such that for every $i\in I:$ $y_{i}^{\ast }\in
A_{i}(x^{\ast })$; $A_{i}(x^{\ast })\cap P_{i}(y^{\ast })=\phi $ and $%
d(x^{\ast },y_{i}^{\ast })=d(X,Y_{i}),$ that is, there exists $(x^{\ast
},y^{\ast })\in X\times Y$ such that for every $i\in I:$

i) $y_{i}^{\ast }\in A_{i}(x^{\ast })$;

ii) sup$_{z\in A_{i}(x^{\ast })}\psi _{i}(y^{\ast },z)\leq 0$;

iii) $d(x^{\ast },y_{i}^{\ast })=d(X,Y_{i}).$ $\square \medskip $

Theorem 11 can be easily proved by using Theorem 10. It refers to the
existence of the solutions for a system of quasi-variational inequalities of
type (2).

\begin{theorem}
\textit{Let }$X$\textit{\ and }$Y_{i}$\textit{\ be non-empty compact and
convex subsets of a normed linear space }$E$\textit{\ for each }$i\in
I=\{1,2,...,n\}.$ Assume that for each $i\in I$, the following conditions
are fulfilled:
\end{theorem}

\textit{(1) The correspondence} $A_{i}:X\rightarrow 2^{Y_{i}}$\textit{\ is
such that }$A_{i}$\textit{\ is lower semicontinuous in }$X^{0},$\textit{\ it
has the property * and }$A_{i}(x)$\textit{\ is a non-empty closed convex
subset of }$Y_{i};$

(2\textit{)\ }$G_{i}:Y_{i}\rightarrow Y^{\prime }$ \textit{is monotone\ with
non-empty values} and $G_{i}:L\cap Y_{i}\rightarrow 2^{E^{\prime }}$\textit{%
\ is lower semicontinuous from the relative topology of }$Y$\textit{\ into
the weak}$^{\ast }-$\textit{topology }$\sigma (E^{\prime },E)$ \textit{of }$%
E^{\prime }$\textit{\ for each one-dimensional flat }$L\subset E.$

\ \textit{Then, there exists }$(x^{\ast },y^{\ast })\in X\times Y$\textit{\
such that }for every $i\in I$:

i) $y_{i}^{\ast }\in A_{i}(x^{\ast });$

ii) $d(x^{\ast },y_{i}^{\ast })=d(X,Y_{i});$

iii) \textit{sup}$_{u\in G_{i}(y_{i}^{\ast })}$Re$\langle u,y_{i}^{\ast
}-z\rangle ]\leq 0$\textit{\ for all }$z\in A_{i}(x^{\ast })$\textit{.}$%
\medskip $

\textit{Proof.} Let us define $\psi _{i}:Y\times Y_{i}\rightarrow R\cup
\{-\infty ,+\infty \}$ by

$\psi _{i}(y,z)=\sup_{u\in G_{i}(z)}$\textit{Re}$\langle u,y_{i}-z\rangle $
for each $(y,z)\in Y\times Y_{i}.$

We have that $y_{i}\notin \{z\in Y_{i}:\psi _{i}(y,z)>0\}$\textit{\ }for
each fixed\textit{\ }$y\in Y$ and, as a consequence of assumption 2), it
follows that for each $y\in Y,$ $\psi _{i}(y,\cdot )$ is concave$.$

All the hypotheses of Theorem 10 are satisfied. According to Theorem 10,
there exists $(x^{\ast },y^{\ast })\in X\times Y$ such that $y_{i}^{\ast
}\in A_{i}(x^{\ast }),$ $d(x^{\ast },y_{i}^{\ast })=d(X,Y_{i})$ for every $%
i\in I$

and

(a) \ \ sup$_{z\in A_{i}(x^{\ast })}\sup_{u\in G_{i}(z)}\mathit{Re}\langle
u,y_{i}^{\ast }-z\rangle \leq 0$ for every $i\in I.$

Finally, we will prove that

sup$_{z\in A_{i}(x^{\ast })}\sup_{u\in G_{i}(y_{i}^{\ast })}\mathit{Re}%
\langle u,y_{i}^{\ast }-z\rangle \leq 0$ for every $i\in I.$

In order to do that, let us consider $i\in I.$

Let $y\in A_{i}(x^{\ast })$, $\lambda \in \lbrack 0,1]$ and $z_{\lambda
}^{i}:=\lambda z+(1-\lambda )y_{i}^{\ast }.$ According to assumption 1), $%
z_{\lambda }^{i}\in A_{i}(x^{\ast }).$

According to (a), we have $\sup_{u\in G_{i}(z_{\lambda }^{i})}\mathit{Re}%
\langle u,y_{i}^{\ast }-z_{\lambda }^{i}\rangle \leq 0$ for each $\lambda
\in \lbrack 0,1]$.

For each $\lambda \in \lbrack 0,1]$, we have that

$t\{\sup_{u\in G_{i}(z_{\lambda }^{i})}\mathit{Re}\langle u,y_{i}^{\ast
}-z\rangle \}=\sup_{u\in G_{i}(z_{\lambda }^{i})}t\mathit{Re}\langle
u,y_{i}^{\ast }-z)\rangle =$

$=\sup_{u\in G_{i}(z_{\lambda }^{i})}\mathit{Re}\langle u,y_{i}^{\ast
}-z_{\lambda }^{i}\rangle \leq 0.$

It follows that for each $\lambda \in \lbrack 0,1],$

(b) $\sup_{u\in G_{i}(z_{\lambda }^{i})}\mathit{Re}\langle u,y_{i}^{\ast
}-z\rangle \leq 0.$

Now, we are using the lower semicontinuity of $G_{i}:L\cap Y_{i}\rightarrow
2^{Y^{\prime }}$ in order to show the conclusion. For each $v_{0}\in
G_{i}(y_{i}^{\ast })$ and $e>0$ let us consider $U_{v_{0}}^{i},$ the
neighborhood of $v_{0}$ in the topology $\sigma (Y^{\prime },Y),$ defined by 
$U_{v_{0}}^{i}:=\{v\in Y^{\prime }:|\func{Re}\langle v_{0}-v,y_{i}^{\ast
}-z\rangle |<e\}.$ As $G_{i}:L\cap Y_{i}\rightarrow 2^{Y^{\prime }}$ is
lower semicontinuous, where $L=\{z_{\lambda }^{i}:\lambda \in \lbrack 0,1]\}$
and $U_{v_{0}}^{i}\cap G_{i}(y_{i}^{\ast })\neq \emptyset ,$ there exists a
non-empty neighborhood $N(y_{i}^{\ast })$ of $y_{i}^{\ast }$ in $L$ such
that for each $z_{i}\in N(y_{i}^{\ast }),$ we have that $U_{v_{0}}^{i}\cap
G_{i}(z_{i})\neq \emptyset .$ Then there exists $\delta \in (0,1],$ $t\in
(0,\delta )$ and $u\in G_{i}(z_{\lambda }^{i})\cap U_{v_{0}}^{i}\neq
\emptyset $ such that $\mathit{Re}\langle v_{0}-u,y_{i}^{\ast }-z\rangle <e.$
Therefore, $\mathit{Re}\langle v_{0},y_{i}^{\ast }-z\rangle <\mathit{Re}%
\langle u,y_{i}^{\ast }-z\rangle +e.$

It follows that $\mathit{Re}\langle v_{0},y_{i}^{\ast }-z\rangle <\mathit{Re}%
\langle u,y_{i}^{\ast }-z\rangle +e<e.$

The last inequality comes from (b). Since $e>0$ and $v_{0}\in
G_{i}(y_{i}^{\ast })$ have been chosen arbitrarily, the next relation holds:

$\mathit{Re}\langle v_{0},y_{i}^{\ast }-z\rangle <0.$

Hence, for each $i\in I,$ we have that $\sup_{u\in G_{i}(y_{i}^{\ast })}%
\mathit{Re}\langle u,y_{i}^{\ast }-z\rangle \leq 0$ for every $z\in
A_{i}(x^{\ast })$.$\square \medskip $

\section{A fuzzy approach to the systems of generalized quasi-variational
inequalities}

\subsection{The free abstract fuzzy economy model}

In this section we extend the results stated in Section 4 in a fuzzy
framework. Firstly, we present the free abstract fuzzy economy model and the
notion of fuzzy equilibrium pair. In the second subsection we will present
the Q'-correspondences and we will prove a theorem of existence of a fuzzy
equilibrium pair for a free abstract economy with Q'-majorized preference
correspondences. We will apply that result to the systems of
quasi-variational inequalities.

We introduced the following model of an abstract fuzzy economy in [17].

Let $I$ be a nonempty set (the set of agents). For each $i\in I$, let $X_{i}$
be a non-empty set of manufacturing commodities, and $Y_{i}$ be a non-empty
set of selling commodities. Define $X:=\underset{i\in I}{\prod }X_{i}$; let $%
A_{i},B_{i}:X\rightarrow \mathcal{F}(Y_{i})$ be the fuzzy constraint
correspondences, $P_{i}:Y:=\underset{i\in I}{\prod }Y_{i}\rightarrow 
\mathcal{F}(Y_{i})$ the fuzzy preference correspondence, $%
a_{i},b_{i}:X\rightarrow (0,1]$ fuzzy constraint functions and $%
p_{i}:Y\rightarrow (0,1]$ fuzzy preference function. We consider that $X_{i}$
and $Y_{i}$ are non-empty subsets of a normed linear space $E$.

\begin{definition}
A \textit{free} \textit{abstract fuzzy economy} is defined as an ordered
family $\Gamma =(X_{i},Y_{i},A_{i},P_{i},a_{i},p_{i})_{i\in I}$.\medskip
\end{definition}

\begin{definition}
A \textit{fuzzy} \textit{equilibrium} pair for $\Gamma $ is defined as a
pair of points $(x^{\ast },y^{\ast })\in X\times Y$ such that for each $i\in
I$, $y_{i}^{\ast }\in (A_{i_{x^{\ast }}})_{a_{i}(x^{\ast })}$ with $%
d(x_{i}^{\ast },y_{i}^{\ast })=d(X_{i},Y_{i})$ and $(A_{i_{x^{\ast
}}})_{a_{i}(x^{\ast })}\cap (P_{i_{x^{\ast }}})_{p_{i}(x^{\ast })}=\emptyset
,$ where $(A_{i_{x^{\ast }}})_{a_{i}(x^{\ast })}=\{z\in Y_{i}:A_{i_{x^{\ast
}}}(z)\geq a_{i}(x^{\ast })\}$ and $(P_{i_{x^{\ast }}})_{p_{i}(x^{\ast
})}=\{z\in Y_{i}:P_{i_{x^{\ast }}}(z)\geq p_{i}(x^{\ast })\}$.
\end{definition}

If $A_{i},P_{i}:X\rightarrow 2^{Y_{i}}$ are classical correspondences, then
we get the definition of free abstract economy and equilibrium pair defined
by W.K. Kim and K. H. Lee in [8].

Whenever $X_{i}=X$ for each $i\in I,$ for the simplicity, we may assume $%
A_{i}:X\rightarrow $ $\mathcal{F}(Y_{i})$ instead of $\ A_{i}:\underset{i\in
I}{\prod }X_{i}\rightarrow \mathcal{F}(Y_{i})$ for the free abstract fuzzy
economy $\Gamma =(X,Y_{i},A_{i},P_{i},a_{i},p_{i})_{i\in I}$ and equilibrium
pair. In particular, when $I=\{1,2,...,n\},$ we may call $\Gamma $ a free
n-person fuzzy game.

The economic interpretation of an \textit{equilibrium} pair for $\Gamma $ is
based on the requirement that for each $i\in I,$ minimize the travelling
cost $d(x_{i},y_{i})$, and also, maximize the preference $P_{i_{y}}$ on the
constraint set $A_{i_{y}}$. Therefore, it is contemplated to find a pair of
points $(x^{\ast },y^{\ast })\in X\times Y$ such that for each $i\in I$, $%
y_{i}^{\ast }\in (A_{i_{x^{\ast }}})_{a_{i}(x^{\ast })}$, $(A_{i_{x^{\ast
}}})_{a_{i}(x^{\ast })}\cap (P_{i_{x^{\ast }}})_{p_{i}(x^{\ast })}=\emptyset 
$ and $\parallel x_{i}^{\ast }-y_{i}^{\ast }\parallel =d(X_{i},Y_{i}),$
where $d(X_{i},Y_{i})=\inf \left\{ \parallel x_{i}^{\ast }-y_{i}^{\ast
}\parallel \mid x_{i}\in X_{i},\text{ }y_{i}\in Y_{i}\right\} $.

When in addition $X_{i}=Y_{i}$ and $A_{i},P_{i}:X\rightarrow 2^{Y_{i}}$ are
classical correspondences for each $i\in I,$ then the previous definitions
can be reduced to the standard definitions of equilibrium theory in
economics due to Yannelis and Prabhakar [24].

\subsection{Fuzzy equilibrium existence and applications to the systems of
generalized quasi-variational inequalities}

In order to prove the theorems in this subsection, we will use the following
results concerning $Q^{\prime }-$majorized correspondences we introduced in
[16].

\begin{definition}
(see [16] Let $X$ be a topological space and $Y$ be a non-empty subset of $\,
$a vector space E, $\theta :X\rightarrow E$ be a mapping and $%
\,T:X\rightarrow 2^{Y}$ be a correspondence.
\end{definition}

(1) $T$ is said to be\textit{\ of class }$Q_{\theta }^{^{\prime }}$ (or $%
Q^{^{\prime }}$) if

\qquad (a) for each $x\in X$, $\theta (x)\notin \overline{T}(x)$ and

\qquad (b) $T$ is lower semicontinuous with open and convex values in $Y$;

(2) A correspondence $T_{x}:X\rightarrow 2^{Y}$ is said to be a $Q_{\theta
}^{^{\prime }}$\textit{-majorant of }$T$\textit{\ at }$x$ \thinspace if
there exists an open neighborhood $N(x)$ of $x$ such that

\qquad (a) For each $z\in N(x)$, $T(z)\subset T_{x}(z)$ and $\theta
(z)\notin \overline{T_{x}}(z)$

\qquad (b) $T_{x}$ is l.s.c. with open and convex values;

(3) $T$ is said to be $Q_{\theta }^{^{\prime }}$\textit{-majorized} if for
each $x\in X$ with $T(x)\neq \emptyset $ there exists a $Q_{\theta
}^{^{\prime }}$-majorant $T_{x}$ of $T$ at $x$.\smallskip

The following Lemma concerning $Q\prime $-majorized correspondences in
needed.\medskip

\begin{theorem}
(see [16] \textit{Let }$X$\textit{\ be a paracompact topological space and }$%
Y$\textit{\ be a non-empty subset of a vector space }$E$\textit{. Let }$%
\theta :X\rightarrow E$\textit{\ be a function and }$T:X\rightarrow
2^{Y}\backslash \{\emptyset \}$\textit{\ be a }$\mathit{Q}^{^{\prime }}$%
\textit{-majorized correspondence. Then, there exists a correspondence }$%
S:X\rightarrow 2^{Y}$\textit{\ of class }$Q^{^{\prime }}$\textit{such that }$%
T(x)\subset S(x)$\textit{\ for each }$x\in X$\textit{\ .}$\medskip $
\end{theorem}

Theorem 13 is an existence theorem of pair equilibrium for a free n person
fuzzy game with upper semi-continuous constraint correspondences and $%
Q_{\theta }^{\prime }$-majorized preference correspondences.

\begin{theorem}
\textit{Let }$\Gamma =(X,Y_{i},A_{i},B_{i},P_{i},a_{i},b_{i},p_{i})_{i\in I}$%
\textit{\ be a free n-person fuzzy game such that for each }$i\in
I=\{1,2...n\}$:
\end{theorem}

(1)\textit{\ }$X$\textit{\ and }$\mathit{Y}_{i}$\textit{\ are non-empty
compact and convex subsets of normed linear space }$\mathit{E}$\textit{;}

(2)\textit{\ }$A_{i}:X\rightarrow $ $\mathcal{F}(Y_{i})$ \textit{is such that%
} $x\rightarrow (A_{i_{x}})_{a_{i}(x)}:X\rightarrow 2^{Y_{i}}$ \textit{is
upper semicontinuous in }$X^{0}$,\textit{\ }$(A_{i_{x}})_{a_{i}(x)}$\textit{%
\ is a non-empty closed convex subset of }$Y_{i}$\textit{\ and }$%
(A_{i_{x}})_{a_{i}(x)}\cap Y_{i}^{0}\neq \emptyset $ for each $x\in X^{0};$

(3)\textit{\ }$P_{i}:Y:\underset{i\in I}{=\prod }Y_{i}\rightarrow \mathcal{F}%
(Y_{i})$ \textit{is such that }$y\rightarrow
(P_{i_{y}})_{p_{i}(y)}:Y\rightarrow 2^{Y_{i}}$ is $Q_{\pi _{i}}^{\prime }-$%
\textit{majorized}$;$

(4) $(P_{i_{y}})_{p_{i}(y)}$ \textit{is nonempty for each} $y\in Y$;\qquad
\qquad \qquad \qquad \qquad \qquad \qquad \qquad \qquad \qquad \qquad \qquad
\qquad \qquad \qquad \qquad \qquad \qquad \qquad \qquad \qquad \qquad \qquad
\qquad \qquad \qquad \qquad \qquad \qquad \qquad \qquad \qquad \qquad \qquad
\qquad \qquad \qquad \qquad \qquad \qquad \qquad \qquad \qquad \qquad \qquad
\qquad \qquad \qquad \qquad \qquad \qquad \qquad \qquad \qquad \qquad \qquad
\qquad \qquad \qquad \qquad \qquad \qquad \qquad \qquad \qquad \qquad \qquad
\qquad \qquad \qquad \qquad \qquad \qquad \qquad \qquad \qquad \qquad \qquad
\qquad \qquad \qquad \qquad \qquad \qquad \qquad \qquad \qquad \qquad \qquad
\qquad \qquad \qquad \qquad \qquad \qquad \qquad \qquad \qquad \qquad \qquad
\qquad \qquad \qquad \qquad \qquad \qquad \qquad \qquad \qquad \qquad \qquad
\qquad \qquad \qquad \qquad \qquad \qquad \qquad \qquad \qquad \qquad \qquad
\qquad \qquad \qquad \qquad \qquad \qquad \qquad \qquad \qquad \qquad \qquad
\qquad \qquad \qquad \qquad \qquad \qquad \qquad \qquad \qquad \qquad \qquad
\qquad \qquad \qquad \qquad \qquad \qquad \qquad \qquad \qquad \qquad \qquad
\qquad \qquad \qquad \qquad \qquad \qquad \qquad \qquad \qquad \qquad \qquad
\qquad \qquad \qquad \qquad \qquad \qquad \qquad \qquad \qquad \qquad \qquad
\qquad \qquad \qquad \qquad \qquad \qquad \qquad \qquad \qquad \qquad \qquad
\qquad \qquad \qquad \qquad \qquad \qquad \qquad \qquad \qquad \qquad \qquad
\qquad \qquad \qquad \qquad \qquad \qquad \qquad \qquad \qquad \qquad \qquad
\qquad \qquad \qquad \qquad \qquad \qquad \qquad \qquad \qquad \qquad \qquad
\qquad \qquad \qquad \qquad \qquad \qquad \qquad \qquad \qquad \qquad \qquad
\qquad \qquad \qquad \qquad \qquad \qquad \qquad \qquad \qquad \qquad \qquad
\qquad \qquad \qquad \qquad \qquad \qquad \qquad \qquad \qquad \qquad \qquad
\qquad \qquad \qquad \qquad \qquad \qquad \qquad \qquad \qquad \qquad \qquad
\qquad \qquad \qquad \qquad \qquad \qquad \qquad \qquad \qquad \qquad \qquad
\qquad \qquad \qquad \qquad \qquad \qquad \qquad \qquad \qquad \qquad \qquad
\qquad \qquad \qquad \qquad \qquad \qquad \qquad \qquad \qquad \qquad \qquad
\qquad \qquad \qquad \qquad \qquad \qquad \qquad \qquad \qquad \qquad \qquad
\qquad \qquad \qquad \qquad \qquad \qquad \qquad \qquad \qquad \qquad \qquad
\qquad \qquad \qquad \qquad \qquad \qquad \qquad \qquad \qquad \qquad \qquad
\qquad \qquad \qquad \qquad \qquad \qquad \qquad \qquad \qquad \qquad \qquad
\qquad \qquad \qquad \qquad \qquad \qquad \qquad \qquad \qquad \qquad \qquad
\qquad \qquad \qquad \qquad \qquad \qquad \qquad \qquad \qquad \qquad \qquad
\qquad \qquad \qquad \qquad \qquad \qquad \qquad \qquad \qquad \qquad \qquad
\qquad \qquad \qquad \qquad \qquad \qquad \qquad \qquad \qquad \qquad \qquad
\qquad \qquad \qquad \qquad \qquad \qquad \qquad \qquad \qquad \qquad \qquad
\qquad \qquad \qquad \qquad \qquad \qquad \qquad \qquad \qquad \qquad \qquad
\qquad \qquad \qquad \qquad \qquad \qquad \qquad \qquad \qquad \qquad \qquad
\qquad \qquad \qquad \qquad \qquad \qquad \qquad \qquad \qquad \qquad \qquad
\qquad \qquad \qquad \qquad \qquad \qquad \qquad \qquad \qquad \qquad \qquad
\qquad \qquad \qquad \qquad \qquad \qquad \qquad \qquad \qquad \qquad \qquad
\qquad \qquad \qquad \qquad \qquad \qquad \qquad \qquad \qquad \qquad \qquad
\qquad \qquad \qquad \qquad \qquad \qquad \qquad \qquad \qquad \qquad \qquad
\qquad \qquad \qquad \qquad \qquad \qquad \qquad \qquad \qquad \qquad \qquad
\qquad \qquad \qquad \qquad \qquad \qquad \qquad \qquad \qquad \qquad \qquad
\qquad \qquad \qquad \qquad \qquad \qquad \qquad \qquad \qquad \qquad \qquad
\qquad \qquad \qquad \qquad \qquad \qquad \qquad \qquad \qquad \qquad \qquad
\qquad \qquad \qquad \qquad \qquad \qquad \qquad \qquad \qquad \qquad \qquad
\qquad \qquad \qquad \qquad \qquad \qquad \qquad \qquad \qquad \qquad \qquad
\qquad \qquad \qquad \qquad \qquad \qquad \qquad \qquad \qquad \qquad \qquad
\qquad \qquad \qquad \qquad \qquad \qquad \qquad \qquad \qquad \qquad \qquad
\qquad \qquad \qquad \qquad \qquad \qquad \qquad \qquad \qquad \qquad \qquad
\qquad \qquad \qquad \qquad \qquad \qquad \qquad \qquad \qquad \qquad \qquad
\qquad \qquad \qquad \qquad \qquad \qquad \qquad \qquad \qquad \qquad \qquad
\qquad \qquad \qquad \qquad \qquad \qquad \qquad \qquad \qquad \qquad \qquad
\qquad \qquad \qquad \qquad \qquad \qquad \qquad \qquad \qquad \qquad \qquad
\qquad \qquad \qquad \qquad \qquad \qquad \qquad \qquad \qquad \qquad \qquad
\qquad \qquad \qquad \qquad \qquad \qquad \qquad \qquad \qquad \qquad \qquad
\qquad \qquad \qquad \qquad \qquad \qquad \qquad \qquad \qquad \qquad \qquad
\qquad \qquad \qquad \qquad \qquad \qquad \qquad \qquad \qquad \qquad \qquad
\qquad \qquad \qquad \qquad \qquad \qquad \qquad \qquad \qquad \qquad \qquad
\qquad \qquad \qquad \qquad \qquad \qquad \qquad \qquad \qquad \qquad \qquad
\qquad \qquad \qquad \qquad \qquad \qquad \qquad \qquad \qquad \qquad \qquad
\qquad \qquad \qquad \qquad \qquad \qquad \qquad \qquad \qquad \qquad \qquad
\qquad \qquad \qquad \qquad \qquad \qquad \qquad \qquad \qquad \qquad \qquad
\qquad \qquad \qquad \qquad \qquad \qquad \qquad \qquad \qquad \qquad \qquad
\qquad \qquad \qquad \qquad \qquad \qquad \qquad \qquad \qquad \qquad \qquad
\qquad \qquad \qquad \qquad \qquad \qquad \qquad \qquad \qquad \qquad \qquad
\qquad \qquad \qquad \qquad \qquad \qquad \qquad \qquad \qquad \qquad \qquad
\qquad \qquad \qquad \qquad \qquad \qquad \qquad \qquad \qquad \qquad \qquad
\qquad \qquad \qquad \qquad \qquad \qquad \qquad \qquad \qquad \qquad \qquad
\qquad \qquad \qquad \qquad \qquad \qquad \qquad \qquad \qquad \qquad \qquad
\qquad \qquad \qquad \qquad \qquad \qquad \qquad \qquad \qquad \qquad \qquad
\qquad \qquad \qquad \qquad \qquad \qquad \qquad \qquad \qquad \qquad \qquad
\qquad \qquad \qquad \qquad \qquad \qquad \qquad \qquad \qquad \qquad \qquad
\qquad \qquad \qquad \qquad \qquad \qquad \qquad \qquad \qquad \qquad \qquad
\qquad \qquad \qquad \qquad \qquad \qquad \qquad \qquad \qquad \qquad \qquad
\qquad \qquad \qquad \qquad \qquad \qquad \qquad \qquad \qquad \qquad \qquad
\qquad \qquad \qquad \qquad \qquad \qquad \qquad \qquad \qquad \qquad \qquad
\qquad \qquad \qquad \qquad \qquad \qquad \qquad \qquad \qquad \qquad \qquad
\qquad \qquad \qquad \qquad \qquad \qquad \qquad \qquad \qquad \qquad \qquad
\qquad \qquad \qquad \qquad \qquad \qquad \qquad \qquad \qquad \qquad \qquad
\qquad \qquad \qquad \qquad \qquad \qquad \qquad \qquad \qquad \qquad \qquad
\qquad \qquad \qquad \qquad \qquad \qquad \qquad \qquad \qquad \qquad \qquad
\qquad \qquad \qquad \qquad \qquad \qquad \qquad \qquad \qquad \qquad \qquad
\qquad \qquad \qquad \qquad \qquad \qquad \qquad \qquad \qquad \qquad \qquad
\qquad \qquad \qquad \qquad \qquad \qquad \qquad \qquad \qquad \qquad \qquad
\qquad \qquad \qquad \qquad \qquad \qquad \qquad

\textit{Then there exists a fuzzy equilibrium pair} \textit{of points} $%
(x^{\ast },y^{\ast })\in X\times Y$ \textit{such that for each }$i\in I$%
\textit{, }$y_{i}^{\ast }\in (A_{i_{x^{\ast }}})_{a_{i}(x^{\ast })}$ \textit{%
with} $d(x_{i}^{\ast },y_{i}^{\ast })=d(X_{i},Y_{i})$ \textit{and} $%
(A_{i_{x^{\ast }}})_{a_{i}(x^{\ast })}\cap (P_{i_{y^{\ast
}}})_{p_{i}(y^{\ast })}=\emptyset .$

\textit{Proof.} Since $x\rightarrow (A_{i_{x}})_{a_{i}(x)}$ satisfies the
whole assumption of Theorem $5$ for each $i\in I$, there exists a point $%
x^{\ast }\in X$ satisfying the system of best proximity pairs, i.e. $\left\{
x^{\ast }\right\} \times (A_{i_{x^{\ast }}})_{b_{i}(x^{\ast })}\subseteq
X\times Y_{i}$ such that $d(x^{\ast },(A_{i_{x^{\ast }}})_{a_{i}(x^{\ast
})})=d(X,Y_{i})$ for each $i\in I.$ Let $\mathcal{A}_{i}:=\left\{ y_{i}\in
(A_{i_{x^{\ast }}})_{a_{i}(x^{\ast })}\text{ / }d(x^{\ast
},y_{i})=d(X,Y_{i})\right\} $ the non-empty best proximity set of the
correspondence $x\rightarrow (A_{i_{x}})_{a_{i}(x)}.$ The set $\mathcal{A}%
_{i}$ is nonempty, closed and convex.

Since $y\rightarrow (P_{i_{y}})_{p_{i}(y)}$ is $Q_{\pi _{i}}^{\prime }$%
\textit{-}majorized for each $i\in I$, according to Theorem 12, we have that
there exists a correspondence $\varphi _{i}:Y\rightarrow 2^{Y_{i}}$ of class 
$Q_{\pi _{i}}^{\prime }$ such that $(P_{i_{y}})_{p_{i}(y)}\subset \varphi
_{i}(y)$ for each $y\in Y.$ Then, $\varphi _{i}$ is lower semicontinuous
with open, convex values and $\pi _{i}(y)\notin \overline{\varphi }_{i}(y)$
for each $y\in Y.$

For each $i\in I$ define a correspondence

$\Phi _{i}:Y\rightarrow 2^{Y_{i}}$ by

$\Phi _{i}(y):=\left\{ 
\begin{array}{c}
\varphi _{i}(y),\text{ \ \ \ \ \ if }y_{i}\notin \mathcal{A}_{i}, \\ 
(A_{i_{x^{\ast }}})_{a_{i}(x^{\ast })}\cap \varphi _{i}(y),\text{ \ if }%
y_{i}\in \mathcal{A}_{i}.%
\end{array}%
\right. $

According to Lemma 1, $\Phi _{i}$ is lower semicontinuous, has convex
values, and $\pi _{i}(y)\notin \overline{\Phi }_{i}(y).$ By applying Theorem 
$3$ to $(Y_{i},\Phi _{i})_{i\in I}$, there exists a maximal element $y^{\ast
}\in Y$ such that $\Phi _{i}\left( y^{\ast }\right) =\emptyset $ for each $%
i\in I.$ For each $y\in Y$ with $y_{i}\notin \mathcal{A}_{i},$ $\Phi _{i}(y)$
is a non-empty subset of $Y_{i}$ because $(P_{i_{y}})_{p_{i}(y)}\neq
\emptyset .$ We have that $y_{i}^{\ast }\in \mathcal{A}_{i}$ and $%
(A_{i_{x^{\ast }}})_{a_{i}(x^{\ast })}\cap \varphi _{i}(y^{\ast })=\emptyset
.$ Since $(P_{i_{y^{\ast }}})_{p_{i}(y^{\ast })}\subset \varphi _{i}(y^{\ast
}),$ it follows that $(A_{i_{x^{\ast }}})_{a_{i}(x^{\ast })}\cap
(P_{i_{y^{\ast }}})_{p_{i}(y^{\ast })}=\emptyset .$ Hence, $y_{i}^{\ast }\in 
\mathcal{A}_{i},$ i.e. $y_{i}^{\ast }\in (A_{i_{x^{\ast }}})_{a_{i}(x^{\ast
})}$ and $d(x^{\ast },y_{i}^{\ast })=d(X,Y_{i})$ for each $i\in I.$ Then $%
(x^{\ast },y^{\ast })$ is a fuzzy equilibrium pair for $\Gamma .$ $\square $

As a consequence of the above theorem, we obtain the next theorem. The
systems of quasi-variational inequalities of type (3) are studied.

\begin{theorem}
\textit{Let }$X$\textit{\ and }$Y_{i}$\textit{\ be non-empty compact and
convex subsets of a normed linear space }$E$\textit{\ for each }$i\in
I=\{1,2,...,n\}.$ Assume that for each $i\in I,$ the following conditions
are fulfilled:
\end{theorem}

\textit{(1) }$A_{i}:X\rightarrow $ $\mathcal{F}(Y_{i})$ \textit{is such that}
$x\rightarrow (A_{i_{x}})_{a_{i}(x)}:X\rightarrow 2^{Y_{i}}$ \textit{is
upper semicontinuous in }$X^{0}$,\textit{\ }$(A_{i_{x}})_{a_{i}(x)}$\textit{%
\ is a non-empty closed convex subset of }$Y_{i}$\textit{, }$%
(A_{i_{x}})_{a_{i}(x)}\cap Y_{i}^{0}\neq \emptyset $ for each $x\in X^{0};$

\textit{The function\ }$\psi _{i}:Y\times Y_{i}\rightarrow \mathbb{R}\cup
\{-\infty ,+\infty \}$\textit{\ is such that }

\ \ \ (\textit{2) for each fixed }$y_{i}\in Y_{i},$ $\psi _{i}(\cdot ,y_{i})$
\textit{is lower semicontinuous such that} $\{z\in Y_{i}:\psi _{i}(y,z)>0\}$ 
\textit{is non-empty for each }$y\in Y$\textit{;}

\ \ \ (\textit{3) }$y_{i}\notin \{z\in Y_{i}:\psi _{i}(y,z)>0\}$\textit{\
for each fixed }$y\in Y;$

\ \ \ (\textit{4) for each} $y\in Y,$ $\psi _{i}(y,\cdot )$ \textit{is
concave.}

\textit{Then, there exists }$(x^{\ast },y^{\ast })\in X\times Y$\textit{\
such that for every }$i\in I$\textit{,}

i) $y_{i}^{\ast }\in (A_{i}(x^{\ast }))_{a_{i}(x^{\ast })};$

ii)$d(x^{\ast },y_{i}^{\ast })=d(X,Y_{i});$

iii) $\sup_{z_{i}^{\ast }\in (A_{i}(x^{\ast }))_{a_{i}(x^{\ast })}}\psi
_{i}(y^{\ast },z_{i}^{\ast })\leq 0.$

\textit{Proof. }For every $i\in I,$ let $P_{i}:Y\rightarrow \mathcal{F}%
(Y_{i})$ such that $(P_{i}(y))_{p_{i}(y)}=\{z\in Y_{i}:\psi _{i}(y,z)>0\}$
for each $y\in Y.$

We shall show that the free abstract economy $%
G=(X,Y_{i},A_{i},P_{i},a_{i},p_{i})_{i\in I}$ satisfies all hypotheses of
Theorem 13.

According to 2), we have that\textit{\ }$y\rightarrow (P_{i}(y))_{p_{i}(y)}$%
\textit{\ }has open lower sections, hence, it is lower semicontinuous with
nonempty compact values and according to 3), $y_{i}\not\in
(P_{i}(y))_{p_{i}(y)}$ for each $y\in Y.$ Assumption 4) implies that $%
y\rightarrow (P_{i}(y))_{p_{i}(y)}$ has convex values. Then, $y\rightarrow
(P_{i}(y))_{p_{i}(y)}$ is Q'-majorozed.

Thus the free abstract economy $G=(X,Y_{i},A_{i},P_{i},a_{i},p_{i})_{i\in I}$
satisfies all hypotheses of Theorem 13. Therefore, there exists $(x^{\ast
},y^{\ast })\in X\times Y$ such that for every $i\in I:$

$y_{i}^{\ast }\in (A_{i}(x^{\ast }))_{a_{i}(x^{\ast })}$; $(A_{i}(x^{\ast
}))_{a_{i}(x^{\ast })}\cap (P_{i}(y^{\ast }))_{p_{i}(y^{\ast })}=\phi $ and $%
d(x^{\ast },y_{i}^{\ast })=d(X,Y_{i}).$

that is, there exists $(x^{\ast },y^{\ast })\in X\times Y$ such that for
every $i\in I:$

i) $y_{i}^{\ast }\in (A_{i}(x^{\ast }))_{a_{i}(x^{\ast })}$;

ii) sup$_{z\in (A_{i}(x^{\ast }))_{a_{i}(x^{\ast })}}\psi _{i}(y^{\ast
},z)\leq 0$;

iii) $d(x^{\ast },y_{i}^{\ast })=d(X,Y_{i}).$ $\square \medskip $

\section{Systems of random quasi-variational inequalities}

In this section we will study the systems of random quasi-variational
inequalities.

In order to prove the existence theorem of random equilibrium pairs for a
random free abstract economy, we need the following result.

\begin{theorem}
(Leese [10], Corollary, pag 408-409). Let $(\Omega ,\mathcal{F})$ be a
measurable space, $\mathcal{F}$ a Suslin family and $X$ a Suslin space.
Suppose that $A:\Omega \rightarrow 2^{X}$ has non-empty values such that Gr$%
A\in \mathcal{F\otimes B}(X).$ Then, there exists a sequence $%
\{g_{n}\}_{n=1}^{\infty }$ of measurable selections of $A$ such that for
each $\omega \in \Omega ,$ $\{g_{n}(\omega ):n\in N\}$ is dense in $A(\omega
).$\medskip 
\end{theorem}

We introduce the model of the random free abstract economy and we define the
random equilibrium pair.

Let $(\Omega ,\mathcal{F})$ be a measurable space, $\mathit{let}$ $I$ be an
index set. For each $i\in I$, let $X_{i}$ be a non-empty set of
manufacturing commodities, and $Y_{i}$ be a non-empty set of selling
commodities. Define $X:=\underset{i\in I}{\prod }X_{i}$; let $A_{i}:\Omega
\times X\rightarrow 2^{Y_{i}}$ be the random constraint correspondence and $%
P_{i}:\Omega \times Y:=\underset{i\in I}{\Omega \times \prod }%
Y_{i}\rightarrow 2^{Y_{i}}$ the random preference correspondence. We
consider that $X_{i}$ and $Y_{i}$ are non-empty subsets of a normed linear
space $E$.

\begin{definition}
\textit{\ }A \textit{free} \textit{abstract economy} is the family $\Gamma
=((\Omega ,\mathcal{F}),X,Y_{i},A_{i},P_{i})_{i\in I}$.\medskip
\end{definition}

\begin{definition}
\textit{\ }A random \textit{equilibrium} pair for $\Gamma $ is defined as a
pair of \textit{measurable functions} $\varphi ^{1}:\Omega \rightarrow X$
and $\varphi ^{2}:\Omega \rightarrow Y$ \textit{such that for each }$i\in I$%
\textit{, }$\pi _{i}(\varphi ^{2})\in A_{i}(\omega ,\varphi ^{1}(\omega ))$%
\textit{\ with} $d(\varphi ^{1}(\omega ),\pi _{i}(\varphi ^{2}(\omega
)))=d(X,Y_{i})$ \textit{and} $A_{i}(\varphi ^{1}(\omega ))\cap P_{i}(\varphi
^{2}(\omega ))=\emptyset $ for all $\omega \in \Omega $.\textit{\medskip }
\end{definition}

Now, we are proving the existence of the random equilibrium for a random
free abstract economy.

\begin{theorem}
\textit{Let }$\Gamma =\{(\Omega ,\mathcal{F}),X,Y_{i},A_{i},P_{i})_{i\in
I}\} $\textit{\ be a free n-person game such that for each }$i\in
I=\{1,2,...,n\}$:
\end{theorem}

(1)\textit{\ }$X$\textit{\ and }$\mathit{Y}_{i}$\textit{\ are non-empty
compact and convex subsets of normed linear space }$\mathit{E}$\textit{;}

(2$)$ $M_{i}=\{(\omega ,x,y):A_{i}(\omega ,x)\cap P_{i}(\omega ,y)\neq
\emptyset \}\in \mathcal{F\otimes }\mathcal{B}(X)\mathcal{\otimes B}(Y),$ $%
N=\{(\omega ,x,y):d(x(\omega ),y_{i}(\omega ))=d(X,Y_{i})$ for each $i\in
I\}\in \mathcal{F\otimes }\mathcal{B}(X)\mathcal{\otimes B}(Y),$ Grcl$%
A_{i}\in \mathcal{F\otimes }\mathcal{B}(X\times Y_{i});$

(3)\textit{\ For each }$\omega \in \Omega ,$ $A_{i}(\omega ,\cdot
):X\rightarrow 2^{Y_{i}}$ is such that $A_{i}(\omega ,\cdot )$ \textit{is
upper semicontinuous in }$X^{0}$\textit{, }$A_{i}(\omega ,x)$\textit{\ is a
closed convex subset of }$Y_{i}$\textit{\ and }$A_{i}(\omega ,x)\cap
Y_{i}^{0}\neq \emptyset $ for each $x\in X^{0}$;

(4)\textit{\ For each }$\omega \in \Omega ,$ $P_{i}(\omega ,\cdot ):Y:%
\underset{i\in I}{=\prod }Y_{i}\rightarrow 2^{Y_{i}}$ \textit{is }$Q_{\pi
_{i}}-$\textit{majorized}$;$

(5) $P_{i}(\omega ,y)$ \textit{is nonempty for each} $(\omega ,y)\in \Omega
\times Y$.

\textit{Then, there exists the measurable functions} $\varphi ^{1}:\Omega
\rightarrow X$ and $\varphi ^{2}:\Omega \rightarrow Y$ \textit{such that for
each }$i\in I$\textit{, }$\pi _{i}(\varphi ^{2})\in A_{i}(\omega ,\varphi
^{1}(\omega ))$\textit{\ with} $d(\varphi ^{1}(\omega ),\pi _{i}(\varphi
^{2}(\omega )))=d(X,Y_{i})$ \textit{and} $A_{i}(\varphi ^{1}(\omega ))\cap
P_{i}(\varphi ^{2}(\omega ))=\emptyset $ for all $\omega \in \Omega $.%
\textit{\medskip }

\textit{Proof.} For each $i\in I,$ define $\phi _{i}:\Omega \rightarrow
2^{X\times Y}$ by $\phi _{i}(\omega )=\{(x,y)\in X\times Y:A_{i}(x)\cap
P_{i}(y)=\emptyset ,$ $y_{i}\in A_{i}(x)$ and $d(x,y_{i})=d(X,Y_{i})\}$ for
each $\omega \in \Omega .$

We also define $\phi :\Omega \rightarrow 2^{X\times Y}$ by $\phi (\omega
)=\cap _{i\in I}\phi _{i}(\omega )$ for each $\omega \in \Omega .$

Then, by Theorem 7, $\phi (\omega )\neq \emptyset .$

Gr$\phi =((\Omega \times X\times Y\backslash (\cup _{i\in I}M_{i}))\cap $Gr$%
\tprod\nolimits_{i\in I}$cl$A_{i}\cap N$ $\in \mathcal{F}\otimes \mathcal{B}%
(X)\otimes \mathcal{B}(Y).$

It follows that $\phi $ satisfies all the conditions of Theorem 14. By
Theorem 14, there exists a measurable selection $\varphi ^{\prime }:$ $%
\Omega \rightarrow X\times Y$ of $\phi .$ Then, there exists $\varphi
^{1}:\Omega \rightarrow X$ and $\varphi ^{2}:\Omega \rightarrow Y$ such that 
$\varphi ^{\prime }(\omega )=(\varphi ^{1}(\omega ),\varphi ^{2}(\omega ))$
for all $\omega \in \Omega .$

We claim that $\varphi ^{1}$ and $\varphi ^{2}$ are measurable. Let $D$ be a
closed subset of $X.$ Then, $D\times Y$ is a closed subset of $X\times Y.$
As $(\varphi ^{1})^{-1}(D)=\{\omega \in \Omega :\varphi ^{1}(\omega )\in
D\}=\{\omega \in \Omega :\varphi (\omega )\in D\times Y\}\in \mathcal{F},$
it follows that $\varphi ^{1}$ is also measurable. We can prove, in the same
way, that $\varphi ^{2}$ is measurable.

Moreover, we have for each $i\in I$, $\pi _{i}(\varphi ^{2})\in A_{i}(\omega
,\varphi ^{1}(\omega ))$\ with $d(\varphi ^{1}(\omega ),\pi _{i}(\varphi
^{2}(\omega )))=d(X,Y_{i})$ and $A_{i}(\varphi ^{1}(\omega ))\cap
P_{i}(\varphi ^{1}(\omega ))=\emptyset $ for all $\omega \in \Omega $. $%
\square $

As an application of Theorem 16, we state the following result concerning
the systems of random generalized quasi-variational inequalities.

\begin{theorem}
\textit{Let }$X$\textit{\ and }$Y_{i}$\textit{\ be non-empty compact and
convex subsets of a normed linear space }$E$\textit{\ for each }$i\in
I=\{1,2,...,n\}.$ Assume that for each $i\in I,$ the following conditions
are fulfilled:
\end{theorem}

\textit{(1) }$M_{i}=\{(\omega ,x,y):A_{i}(\omega ,x)\cap \{z\in Y_{i}:\psi
_{i}(\omega ,y,z)>0\}\neq \emptyset \}\in F\otimes B(X)\otimes B(Y),$\textit{%
\ }$N=\{(\omega ,x,y):d(x(\omega ),y_{i}(\omega ))=d(X,Y_{i})$\textit{\ for
each }$i\in I\}\in F\otimes B(X)\otimes B(Y),$\textit{\ }Grcl$A_{i}\in
F\otimes B(X\times Y_{i});$

\textit{(2) For each }$\omega \in \Omega ,$\textit{\ }$A_{i}(\omega ,\cdot
):X\rightarrow 2^{Y_{i}}$\textit{\ is such that }$A_{i}(\omega ,\cdot )$%
\textit{\ is upper semicontinuous in }$X^{0},$\textit{\ }$A_{i}(\omega ,x)$%
\textit{\ is a closed convex subset of }$Y_{i}$\textit{\ and }$A_{i}(\omega
,x)\cap Y_{i}^{0}\neq \emptyset $\textit{\ for each }$x\in X^{0}$\textit{;}

\textit{(3) The function\ }$\psi _{i}:\Omega \times Y\times Y_{i}\rightarrow 
\mathbb{R}\cup \{-\infty ,+\infty \}$\textit{\ is such that }

\textit{\ \ \ (a) For each }$\omega \in \Omega $\textit{\ and for each fixed 
}$y_{i}\in Y_{i},$\textit{\ }$\psi _{i}(\omega ,\cdot ,y_{i})$\textit{\ is
lower semicontinuous such that }$\{z\in Y_{i}:\psi _{i}(\omega ,y,z)>0\}$%
\textit{\ is non-empty for each }$\omega \in \Omega $\textit{\ and }$y\in Y$%
\textit{;}

\textit{\ \ \ (b) }$y_{i}\notin \{z\in Y_{i}:\psi _{i}(\omega ,y,z)>0\}$%
\textit{\ for each }$\omega \in \Omega $\textit{\ and for each fixed }$y\in
Y;$

\textit{\ \ \ (c) for each }$\omega \in \Omega $\textit{\ and for each }$%
y\in Y,$\textit{\ }$\psi _{i}(\omega ,y,\cdot )$\textit{\ is concave.}

\textit{Then, there exists the measurable functions }$\varphi ^{1}:\Omega
\rightarrow X$\textit{\ and }$\varphi ^{2}:\Omega \rightarrow Y$\textit{\
such that for each }$i\in I$\textit{, }$\pi _{i}(\varphi ^{2})\in
A_{i}(\omega ,\varphi ^{1}(\omega ))$\textit{\ with }$d(\varphi ^{1}(\omega
),\pi _{i}(\varphi ^{2}(\omega )))=d(X,Y_{i})$\textit{\ and }$\sup_{z_{i}\in
A_{i}(\varphi ^{1}(\omega ))}\psi _{i}(\varphi ^{2}(\omega ),z_{i})\leq 0$%
\textit{\ for all }$\omega \in \Omega $\textit{.\medskip }

\textit{Proof. }Let us fix $\omega \in \Omega .$

For every $i\in I,$ let $P_{i}:Y\rightarrow Y_{i}$ be defined by $%
P_{i}(\omega ,y)=\{z\in Y_{i}:\psi _{i}(\omega ,y,z)>0\}$ for each $y\in Y.$

We shall show that the free abstract economy $G=\{(\Omega ,\mathcal{F}%
),X,Y_{i},A_{i},P_{i})_{i\in I}\}$ satisfies all hypotheses of Theorem 15.

According to 3 a), we have that\textit{\ }$P_{i}(\omega ,\cdot )$\textit{\ }%
is lower semicontinuous and then, Q-majorized with nonempty compact values
and according to 3 b), $y_{i}\not\in P_{i}(\omega ,y)$ for each $y\in Y.$
Assumption 3 c) implies that $P_{i}(\omega ,\cdot )$ has convex values.

Thus the free abstract economy $G=\{(\Omega ,\mathcal{F}%
),X,Y_{i},A_{i},P_{i})_{i\in I}\}$ satisfies all hypotheses of Theorem 15.
Therefore, there exists the measurable functions $\varphi ^{1}:\Omega
\rightarrow X$ and $\varphi ^{2}:\Omega \rightarrow Y$ such that for each $%
i\in I$, $\pi _{i}(\varphi ^{2})\in A_{i}(\omega ,\varphi ^{1}(\omega ))$\
with $d(\varphi ^{1}(\omega ),\pi _{i}(\varphi ^{2}(\omega )))=d(X,Y_{i})$
and $A_{i}(\varphi ^{1}(\omega ))\cap P_{i}(\varphi ^{2}(\omega ))=\emptyset 
$ for all $\omega \in \Omega ,$ that is, there exists the measurable
functions $\varphi ^{1}:\Omega \rightarrow X$ and $\varphi ^{2}:\Omega
\rightarrow Y$ such that for each $i\in I$, $\pi _{i}(\varphi ^{2})\in
A_{i}(\omega ,\varphi ^{1}(\omega ))$\ with $d(\varphi ^{1}(\omega ),\pi
_{i}(\varphi ^{2}(\omega )))=d(X,Y_{i})$ and $\sup_{z_{i}\in A_{i}(\varphi
^{1}(\omega ))}\psi _{i}(\varphi ^{2}(\omega ),z_{i})\leq 0$ for all $\omega
\in \Omega $\textit{. }$\square \medskip $

If $|I|=1$, we obtain the following corollary.

\begin{corollary}
\textit{Let }$X$\textit{\ and }$Y$\textit{\ be non-empty compact and convex
subsets of a normed linear space }$E$ \textit{and assume that the following
conditions are fulfilled:\ }
\end{corollary}

\textit{(1) }$M=\{(\omega ,x,y):A(\omega ,x)\cap \{z\in Y:\psi (\omega
,y,z)>0\}\neq \emptyset \}\in F\otimes B(X)\otimes B(Y),$\textit{\ }$%
N=\{(\omega ,x,y):d(x(\omega ),y(\omega ))=d(X,Y)\}\in F\otimes B(X)\otimes
B(Y),$\textit{\ }Grcl$A\in F\otimes B(X\times Y);$

\textit{(2) For each }$\omega \in \Omega ,$\textit{\ }$A(\omega ,\cdot
):X\rightarrow 2^{Y}$\textit{\ is such that }$A(\omega ,\cdot )$\textit{\ is
upper semicontinuous in }$X^{0},$\textit{\ }$A(\omega ,x)$\textit{\ is a
closed convex subset of }$Y$\textit{\ and }$A(\omega ,x)\cap Y^{0}\neq
\emptyset $\textit{\ for each }$x\in X^{0}$\textit{;}

\textit{(3) The function\ }$\psi :\Omega \times Y\times Y\rightarrow R\cup
\{-\infty ,+\infty \}$\textit{\ is such that }

\textit{\ \ \ (a) For each }$\omega \in \Omega $\textit{\ and for each fixed 
}$y\in Y,$\textit{\ }$\psi (\omega ,\cdot ,y)$\textit{\ is lower
semicontinuous such that }$\{z\in Y:\psi (\omega ,y,z)>0\}$\textit{\ is
non-empty for each }$\omega \in \Omega $\textit{\ and }$y\in Y$\textit{;}

\textit{\ \ \ (b) }$y\notin \{z\in Y:\psi (\omega ,y,z)>0\}$\textit{\ for
each }$\omega \in \Omega $\textit{\ and for each fixed }$y\in Y;$

\textit{\ \ \ (c) for each }$\omega \in \Omega $\textit{\ and for each }$%
y\in Y,$\textit{\ }$\psi (\omega ,y,\cdot )$\textit{\ is concave.}

\textit{Then, there exists the measurable functions }$\varphi ^{1}:\Omega
\rightarrow X$\textit{\ and }$\varphi ^{2}:\Omega \rightarrow Y$\textit{\
such that: }$\pi (\varphi ^{2})\in A(\omega ,\varphi ^{1}(\omega ))$\textit{%
\ with }$d(\varphi ^{1}(\omega ),\pi (\varphi ^{2}(\omega )))=d(X,Y)$\textit{%
\ and }$\sup_{z\in A(\varphi ^{1}(\omega ))}\psi (\varphi ^{2}(\omega
),z)\leq 0$\textit{\ for all }$\omega \in \Omega $\textit{.\medskip }

If $\psi _{i}=0$ in Theorem 17$,$ we obtain the following random fixed point
theorem.

\begin{theorem}
\textit{Let }$X$\textit{\ and }$Y_{i}$\textit{\ be non-empty compact and
convex subsets of a normed linear space }$E$\textit{\ for each }$i\in
I=\{1,2,...,n\}.$ Assume that for each $i\in I,$ the following conditions
are fulfilled:
\end{theorem}

\textit{(1) }$N=\{(\omega ,x,y):d(x(\omega ),y_{i}(\omega ))=d(X,Y_{i})$%
\textit{\ for each }$i\in I\}\in F\otimes B(X)\otimes B(Y),$\textit{\ }Grcl$%
A_{i}\in F\otimes B(X\times Y_{i});$

\textit{(2) For each }$\omega \in \Omega ,$\textit{\ }$A_{i}(\omega ,\cdot
):X\rightarrow 2^{Y_{i}}$\textit{\ is such that }$A_{i}(\omega ,\cdot )$%
\textit{\ is upper semicontinuous in }$X^{0},$\textit{\ }$A_{i}(\omega ,x)$%
\textit{\ is a closed convex subset of }$Y_{i}$\textit{\ and }$A_{i}(\omega
,x)\cap Y_{i}^{0}\neq \emptyset $\textit{\ for each }$x\in X^{0}$\textit{.}

\textit{Then, there exists the measurable functions }$\varphi ^{1}:\Omega
\rightarrow X$\textit{\ and }$\varphi ^{2}:\Omega \rightarrow Y$\textit{\
such that for each }$i\in I$\textit{, }$\pi _{i}(\varphi ^{2})\in
A_{i}(\omega ,\varphi ^{1}(\omega ))$\textit{\ with }$d(\varphi ^{1}(\omega
),\pi _{i}(\varphi ^{2}(\omega )))=d(X,Y_{i}).$\textit{\ \medskip }


\begin{thebibliography}{99}
\bibitem{1} Bello Cruz J. Y. and Iusem A. N. (2009). A Strongly Convergent
Direct Method for Monotone Variational Inequalities in Hilbert Spaces. 
\textit{Numerical Functional Analysis and Optimization,} \textbf{30}%
(1-2):23-36.

\bibitem{2} Caruso A. O., Khan A. A. and Raciti F. \ (2009). Continuity
Results for a Class of Variational Inequalities with Applications to
Time-Dependent Network Problems. \textit{Numerical Functional Analysis and
Optimization.} 30(11-12):1272-1288.

\bibitem{3} Cegielski A., Gibali A., Reich S. and Zalas R.. An Algorithm for
Solving the Variational Inequality Problem over the Fixed Point Set of a
Quasi-Nonexpansive Operator in Euclidean Space.\textit{Numerical Functional
Analysis and Optimization.} DOI: 10.1080/01630563.2013.771656

\bibitem{4} Cegielski A. and Zalas R.\ \ (2013). Methods for Variational
Inequality Problem Over the Intersection of Fixed Point Sets of
Quasi-Nonexpansive Operators. \textit{Numerical Functional Analysis and
Optimization.} 34(3):255-283.

\bibitem{4} Ding X. P. (1998). Equilibria of noncompact generalized games
with U-majorized preference correspondences. \textit{Appl. Math. Lett.}
11(5): 115-119.

\bibitem{5} Huang X. X. and Yang X. Q. (2010). Levitin--Polyak
Well-Posedness of Vector Variational Inequality Problems with Functional
Constraints. \textit{Numerical Functional Analysis and Optimization}.
31(4):440-459.

\bibitem{6} Kim W. K. (2006). On general best proximity pairs and
equilibrium pairs in free generalized games. \textit{J. Chungceong Math. Soc.%
} 19(1): 19-26.

\bibitem{7} Kim W. K. and Lee K. H. (2006). Existence of best proximity
pairs and equlibrium pairs. \textit{J. Math. An. Appl.} 316:433-446.

\bibitem{8} Lalitha C. S. and Bhatia G. (2009). Well-Posedness for
Variational Inequality Problems with Generalized Monotone Set-Valued Maps, 
\textit{Numerical Functional Analysis and Optimization,} 30(5-6):548-565.

\bibitem{9} Leese S. J. (1974). Multifunctions of Souslin type, \textit{%
Bull. Austral. Math.} \textit{Soc.} 11:395-411.

\bibitem{10} Li J., C. Zhang and X. Ma (2008). On the Metric Projection
Operator and Its Applications to Solving Variational Inequalities in Banach
Spaces. \textit{Numerical Functional Analysis and Optimization, }%
29(3-4):410-418.

\bibitem{11} Liu X., Cai H. (2001). Maximal Elements and Equilibrium of
Abstract Economy. \textit{Appl. Math. Mech.}, 22(\textit{10)}:1225-1230.

\bibitem{12} Maing\'{e} P.-E. \ (2008). Extension of the Hybrid Steepest
Descent Method to a Class of Variational Inequalities and Fixed Point
Problems with Nonself-Mappings. \textit{Numerical Functional Analysis and
Optimization.} 29(7-8):820-834.

\bibitem{14} Michael E. (1956). Continuous selection I.\textit{\ Annals of
Mathematics.\ }63(2):361-382.

\bibitem{15} Noor M. A. and Elsanousi S. A. (1993). Iterative algorithms for
random variational inequalities, \textit{Pan Amer. Math. Journal.} 3:39-50.

\bibitem{16} Patriche M. (2009). A new fixed point theorem and its
applications in equilibrium theory. \textit{Fixed Point Theory}. 1:159-171.

\bibitem{17} Patriche M. (2009). Equilibria for free abstract fuzzy
economies, \textit{An. St. Ovidius University of Constanta}. 17(2):143-154.

\bibitem{8} Reich S. (1978). Approximate selections, best approximations,
fixed points, and invariant sets.\textit{\ J. Math. Anal. Appl.} 62:104-113.

\bibitem{9} Sehgal V. M. and Singh S. P. \ (1988). A generalization to
multifunctions of Fan's best approximation theorem. \textit{Proc. Amer.
Math. Soc. }102:534-537.

\bibitem{10} Shafer W. and Sonnenschein H. (1975). Equilibrium in abstract
economies without ordered preferences.\textit{\ Journal of Mathematical
Economics} 2:345-348.

\bibitem{11} Srinivasan P. S., Veeramani P. (2003). On best approximation
pair theorems and fixed point theorems\textit{.} \textit{Abstr. Appl. Anal.}
2003: 33-47.

\bibitem{12} Srinivasan P. S., Veeramani P. (2004). On existence of
equilibrium pair for constraint generalized games. \textit{Fixed Point
Theory Appl.} 2004: 21-29.

\bibitem{18} Wu X. (1997). A new fixed point theorem and its applications.%
\textit{\ Proc. Amer. Math. Soc.} 125:1779-1783.

\bibitem{19} Yannelis N. C. and Prabhakar N. D. (1983). Existence of maximal
elements and equilibrium in linear topological spaces. \textit{J. Math.
Econom.} 12:233-245.

\bibitem{0} Yaoa Y., Lioub Y.-C. and Shahzad N. (2012). Construction of
Iterative Methods for Variational Inequality and Fixed Point Problems. 
\textit{Numerical Functional Analysis and Optimization}. 33(10):2012.

\bibitem{20} Yuan v (1998). \textit{The Study of Minimax Inequalities and
Applications to Economies and} \textit{Variational Inequalities, }Memoirs of
the American Society. 132.
\end{thebibliography}
\end{document}